\documentclass[11pt]{article}
\usepackage{amsmath,amsfonts,amssymb,mathrsfs,cases,mathdots,colordvi,url}
\usepackage{enumerate}
\usepackage{amssymb,latexsym}
\usepackage{graphicx}
\usepackage{amstext}
\usepackage{amsmath}
\usepackage{amsopn}
\usepackage{amsbsy}
\usepackage{amscd}
\usepackage{amsxtra}
\usepackage{graphics}
\usepackage{booktabs,multirow}
\usepackage{tabularx}
\usepackage{makecell}
\usepackage{bbding}
\usepackage{pifont}
\usepackage{hyperref}
\usepackage{tikz}
\usepackage{pifont}
\usepackage{amsthm}
\usetikzlibrary{shapes.geometric,arrows}
\usetikzlibrary{positioning, shapes.geometric}
\tikzstyle{block} = [rectangle, draw, fill=white!50,
text width=33em, text centered, rounded corners, minimum height=3em]
\tikzstyle{line} = [draw, -latex']
\tikzstyle{arrow} = [thick,->,>=stealth]
\DeclareMathOperator*{\argmin}{arg\,min}

\newtheorem{theorem}{Theorem}[part]
\newtheorem{definition}{Definition}[part]
\newtheorem{proposition}{Proposition}[part]

\newtheorem{remark}{Remark}[part]
\newtheorem{example}{Example}[part]

\numberwithin{equation}{section}

\begin{document}
\title{Combined approach with second-order optimality conditions for bilevel programming problems
\thanks{This paper is dedicated to Professor Roger J-B Wets on the occasion of his 85th birthday.  Authors listed in alphabetical order.}
\thanks{The research of Ye is supported by NSERC. Zhang's work is supported by National Science Foundation of China 11971220, Shenzhen Science and Technology Program (No. RCYX20200714114700072), the Stable Support Plan Program of Shenzhen Natural Science Fund (No. 20200925152128002). }}

\author{Xiaoxiao Ma\thanks{Department of Mathematics and Statistics, University of Victoria, Victoria, B.C., Canada V8W 2Y2. Email: xiaoxiaoma@uvic.ca.}
\and
Wei Yao\thanks{Department of Mathematics, Southern University of Science and Technology, and National Center for Applied Mathematics Shenzhen, Shenzhen, Guangdong, China. Email: yaow@sustech.edu.cn.}
\and
Jane J. Ye\thanks{Department of Mathematics and Statistics, University of Victoria, Victoria, B.C., Canada V8W 2Y2. Email: janeye@uvic.ca.}
\and
Jin Zhang\thanks{Corresponding author. Department of Mathematics and SUSTech International Center for Mathematics, Southern University of Science and Technology, and National Center for Applied Mathematics Shenzhen, Shenzhen, Guangdong, China. Email: zhangj9@sustech.edu.cn.}}
\date{}

\maketitle

\begin{abstract}
\noindent
In this paper, we propose a combined approach with second-order optimality conditions of the lower level problem to study constraint qualifications and optimality conditions for {bilevel programming problems}. The new method is inspired by the combined approach developed by Ye and Zhu in 2010, where the authors combined the classical first-order and the value function approaches to derive new necessary optimality conditions. In our approach, we add a second-order optimality condition to the combined program as a new constraint.
  We show that when all known approaches fail, adding the second-order optimality condition as a constraint makes the corresponding partial calmness condition  and the resulting necessary optimality condition easier to hold. We also give some discussions on advantages and disadvantages of the combined approaches with the first-order and the second-order information.
\vspace{3mm}

\noindent {\bf Key words:}\hspace{2mm} partial calmness, bilevel program,
  optimality condition, second-order optimality condition \vspace{3mm}

\noindent {\bf AMS Subject Classifications:}\hspace{2mm} 90C26, 90C30, 90C31, 90C33, 90C46, 49J52, 91A65

\end{abstract}

\section{Introduction}
\label{intro}

In this paper we consider the following bilevel programming problem (BLPP):
\begin{equation}\label{BLPP}\tag{BLPP}
\begin{aligned}
\min_{x,y}&\  F(x,y)\\
\mathrm{s.t.}&\ y\in S(x), \  G(x,y)\leq0,
\end{aligned}
\end{equation}
where $S(x)$ denotes the solution set of the lower level program
\begin{equation}\label{Lx}\tag*{$P(x)$}
\min _{y} f(x, y) \quad {\rm s.t. } \quad g(x, y) \leq 0.
\end{equation}
For convenience, we denote the feasible set of \ref{Lx} by
\begin{eqnarray*}
	Y(x):=\left\{y \in \mathbb{R}^{m}: g(x, y) \leq 0\right\}.
\end{eqnarray*}
Here $x \in \mathbb{R}^{n}$, $y \in \mathbb{R}^{m}$ and the mappings $F, f: \mathbb{R}^{n} \times \mathbb{R}^{m} \rightarrow \mathbb{R}$, $G: \mathbb{R}^{n} \times \mathbb{R}^{m} \rightarrow \mathbb{R}^{q}$, $g: \mathbb{R}^{n} \times \mathbb{R}^{m} \rightarrow \mathbb{R}^{p}$.
Unless otherwise specified, we assume that $F, G$ are continuously differentiable and $f,g$ are three times continuously differentiable.

The bilevel programming problem has many applications including the principal-agent moral hazard problem \cite{M7599}, hyperparameter optimization and meta-learning in machine learning \cite{FFSGP18,KBHP08,LMYZZ20,YYZZ}. More applications can be found in \cite{B98,D18,DKPK15,SIB97}. For a comprehensive review, we refer to  \cite{Dempebook} and the references therein.

It is well known that optimality conditions of the lower level program are very useful in the reformulation of BLPPs both theoretically and computationally. {The classical} Karush-Kuhn-Tucker (KKT) approach is to replace the lower level program by its KKT condition and minimize over the original variables as well as multipliers.
When the lower level program has inequality constraints, the KKT reformulation has been studied in {the framework of the mathematical program} with equilibrium constraints (MPEC).  However there are some problems associated with the KKT approach.  First, the KKT condition may  only be necessary but not sufficient.  In this case the KKT approach will enlarge the feasible region  and hence the resulting single level problem may not be equivalent to the original bilevel program.
Moreover an example in  Mirrlees  \cite{M7599} shows that  the solution set of the equivalent single level reformulation  may not even include solutions of the original bilevel program. Second, even in the case where the KKT conditions are necessary and sufficient for $y\in S(x)$,  treating multipliers of the lower level as extra variables can still make the resulting single level reformulation different from the original BLPPs  in the  sense of local optimality; see  \cite{DD12}.  Recently \cite{BM21} has discussed the issue of equivalence for more general problems for which some reformulations may include implicit variables with the BLPP as an example.  Recently  reformulations using the  Bouligand (B-) stationary condition for the lower level program:
 $$0\in \nabla_y f(x,y)+\widehat{N}_{Y(x)}(y),$$
 where $\widehat{N}_C$ denotes the regular normal cone to set $C$ to replace the lower level program
 have been investigated; see \cite{AHO18,GY17,GY19,{KYYZ}}.  Calmness properties for the KKT reformulation and the B-stationarity reformulation have been compared in \cite{AHO18,GY17} in the context of  MPECs and it was discovered that usually the  B-stationarity reformulation is easier to satisfy the calmness condition than the KKT reformulation. Note that extra assumptions (at least the smoothness of the objective function of the lower level program) are always required for a reformulation using optimality conditions for the lower level program.

Contrast to any reformulation using optimality conditions for the lower level program,  the value function approach  proposed by Outrata \cite{O90} for numerical purpose and used by Ye and Zhu \cite{YZ95} for optimality conditions does not require any extra assumptions. By this approach, one defines the value function as an extended real-valued function
\begin{equation*}
V(x):=\inf_y\big\{f(x,y): g(x,y)\leq0\big\},
\end{equation*}
and replaces the original BLPP by the following
equivalent problem:
\begin{equation}\label{VP}\tag{VP}
\begin{aligned}
\min_{x,y}\ & F(x,y) \quad
\mathrm{s.t.}
&f(x,y)-V(x)\leq0,\
g(x,y)\leq0,\
G(x,y)\leq0.
\end{aligned}
\end{equation}
However, since the value function constraint  $f(x,y)-V(x)\leq0$ is actually an equality constraint, the nonsmooth Mangasarian-Fromovitz constraint qualification (MFCQ) for \eqref{VP} will never hold \cite[Proposition 3.2]{YZ95}. To derive necessary optimality conditions for BLPPs, Ye and Zhu \cite[Definition 3.1 and Proposition 3.3]{YZ95} proposed the partial calmness condition for \eqref{VP} under which the difficult constraint $f(x,y)-V(x)\leq0$ was added as a penalty term to the objective function.

Although it was proved in \cite{YZ95} that the partial calmness condition for \eqref{VP} holds automatically for the minmax problem and the bilevel program where the lower level program is linear in both upper and lower level variables, the partial calmness condition for \eqref{VP} is celebrated but has been shown to be restrictive {(cf. \cite{DZ13,HS11,MMZ20,M18,Ye20})}. To improve the value function approach, Ye and Zhu \cite{YZ10} proposed a combination of the classical KKT and the value function approach. The resulting problem is  the combined problem using KKT condition:
\begin{equation}\label{CP}\tag{CP}
\begin{aligned}
\min_{x,y,u}\ & F(x,y) \\
\mathrm{s.t.\ \,} &f(x,y)-V(x)\leq0, \ \nabla_y f(x,y)+\nabla_y g(x,y)^Tu=0, \\
&g(x,y)\leq0,\ u\geq0,\ u^T g(x,y)=0,
\ G(x,y)\leq0.
\end{aligned}
\end{equation}
Problem \eqref{CP} is equivalent to problem \eqref{BLPP} (in the sense of \cite[Proposition 3.1]{YZ10}) when the KKT condition holds at each optimal solution of the lower level program. Similar to \cite{YZ95},  to deal with the fact that  the nonsmooth MFCQ also fails for \eqref{CP}, the corresponding  partial calmness condition for \eqref{CP} was proposed in \cite[Definition 3.1]{YZ10}.

Note that the reformulation \eqref{CP} requires the validity of the KKT conditions at each optimal solution of the lower level program. To deal with the case where the KKT condition may not hold at all the solutions of the lower level program, the Fritz John (FJ) condition was considered. In \cite{AS13}, Allende and Still replaced the lower-level program of BLPPs with the FJ condition (without a value function constraint). Regarding the combined approach, Ke et al. \cite{KYYZ} proposed the following combined program using the FJ condition:
\begin{equation}\label{CPFJ}\tag{CPFJ}
\begin{aligned}
\min_{x,y,u_0, u}\ & F(x,y) \\
\mathrm{s.t.\ \,} &f(x,y)-V(x)\leq0, \ u_0\nabla_y f(x,y)+\nabla_y g(x,y)^Tu=0, \\
&g(x,y)\leq0,\  (u_0,u) \geq0,\ u^T g(x,y)=0, \ \sum_{i=0}^p u_i=1,
\ G(x,y)\leq 0.
\end{aligned}
\end{equation}

The equivalence of problem \eqref{CPFJ} and problem \eqref{BLPP} (for both local and global optimal points) is a special case of \cite[Theorem 4.5]{BM21}.

Similar to Ye and Zhu's work in \cite{YZ10}, Ke et al. proposed the following partial calmness condition for \eqref{CPFJ} in \cite{KYYZ}.
\begin{definition}[Partial calmness for \eqref{CPFJ}]\label{Defn1.2}
	Let $(\bar{x},\bar{y},\bar u_0,\bar{u})$ be a local solution of \eqref{CPFJ}. We say that \eqref{CPFJ} is partially calm at $(\bar{x},\bar{y},\bar u_0,\bar{u})$ if there exists $\mu\geq 0$ such that $(\bar{x},\bar{y},\bar u_0,\bar{u})$ is a local solution of the partially penalized problem:
	\begin{equation}\label{CP2}\tag{\text{CPFJ$_\mu$}}
	\begin{aligned}
	\min_{x,y,u_0, u}\ & F(x,y)+\mu\big(f(x,y)-V(x)\big) \\
	\mathrm{s.t.\ \,}
	&u_0\nabla_y f(x,y)+\nabla_y g(x,y)^Tu=0, \\
	&g(x,y)\leq0,\ (u_0,u)\geq0,\ u^T g(x,y)=0,\ \sum_{i=0}^p u_i=1,
	\ G(x,y)\leq0.
	\end{aligned}
	\end{equation} 	
\end{definition}
Moreover, they analyzed the partial calmness for the combined program based on FJ conditions from a generic point of view and proved that the partial calmness for \eqref{CPFJ} is generic when the upper level variable has dimension one.

The following two combined problems and the corresponding partial calmness conditions are also discussed in \cite{KYYZ}.
\begin{equation}\label{CP_{FJ}}\tag{CP$_{FJ}$}
\begin{aligned}
\min _{x, y}\ & F(x, y) \quad
\text { s.t. }  f(x, y)-V(x) \leq 0, \ (x, y) \in \Sigma_{F J},\ G(x, y) \leq 0,
\end{aligned}
\end{equation}
where $\Sigma_{\mathrm{F J}}:=\big\{(x, y) \in \mathbb{R}^{n+m}:y\text{  satisfies the FJ condition  for }P(x) \big\}$, and
the combined program with the B-stationary condition for the lower level program:
\begin{equation}\label{CPB}\tag{CPB}
\begin{aligned}
\min _{x, y}\ & F(x, y) \quad
\text { s.t. }  f(x, y)-V(x) \leq 0, \ 0 \in \nabla_y f(x, y)+\widehat{N}_{Y(x)}(y), G(x, y) \leq 0.
\end{aligned}
\end{equation}

{Although the partial calmness for the combined program may hold quite often, there are still cases where it does not hold; see e.g. Examples \ref{unc}, \ref{cons1}, \ref{cons2} in this paper}.
The main goal of this paper is to investigate the following question:
\begin{equation}\label{q}\tag{Q}
\begin{aligned}
&\text{How to derive necessary optimality conditions for bilevel problems}\\
&\qquad\text{ when necessary optimality conditions for problems \eqref{VP} and \eqref{CP}  do not hold? }
\end{aligned}
\end{equation}

\textbf{Contributions.}
To answer (Q), we propose to use second-order optimality conditions of the lower level program {on top of the value function constraint and the first-order optimality constraint. The key point of adding the second order condition is to  increase the freedom of choosing multipliers. Since  each extra redundant constraint is  associated with a multiplier, the more redundant constraints we add, the weaker is the optimality condition and hence easier for the resulting necessary optimality condition to hold. In this sense,
the resulting necessary
optimality condition by adding the second-order condition is much more likely to hold than the one adding the first-order condition only, i.e., \eqref{CP},  which is in term more likely to hold than the one without adding any optimality condition, i.e., \eqref{VP}.
}

To illustrate our approach, consider the following KKT combined program:
\begin{equation}\label{KKTCP}\tag{KKTCP}
\begin{aligned}
\min_{x,y} &\  F(x,y) \quad
\mathrm{s.t.} &\ f(x, y)-V(x) \leq 0, \ (x,y)\in \Sigma_{\mathrm{KKT}},\ G(x, y) \leq 0,
\end{aligned}
\end{equation}
where
{ $$\Sigma_{\mathrm{KKT}}:=\left\{(x, y) :\exists u \mbox{ s.t. } \begin{array}{l} \\\nabla_y f(x,y)+\nabla_y g(x,y)^Tu=0, \\
g(x,y)\leq0,\ u\geq0,\ u^T g(x,y)=0\end{array} \right\}$$}
and its partially penalized problem:
\begin{equation}\label{KKTCP-mu}\tag{\text{KKTCP$_\mu$}}
	\begin{aligned}
	\min_{x,y}\ & F(x,y)+\mu\big(f(x,y)-V(x)\big) \quad
	\mathrm{s.t.\ \,}
	&\ (x,y)\in \Sigma_{\mathrm{KKT}},\ G(x, y) \leq 0.
	\end{aligned}	
	\end{equation} 	
	Note that the combined program \eqref{CP} is a relaxed problem of \eqref{KKTCP}  in the sense that the minimization is also performed on multipliers in problem \eqref{CP}.	
	To use the second-order information, we propose the following second-order combined problem:
\begin{equation}\label{SOCP}\tag{SOCP}
\begin{aligned}
\min_{x,y}&\  F(x,y) \quad
\mathrm{s.t.}&\ f(x, y)-V(x) \leq 0, \ (x,y)\in \Sigma_{\mathrm{SOC}},\ G(x, y) \leq 0,
\end{aligned}
\end{equation}
where $$\Sigma_{\mathrm{SOC}}:=\Big\{(x, y) \in \mathbb{R}^{n+m}:y\text{ satisfies a second-order optimality condition for }P(x) \Big\},$$
and its partially penalized problem:
\begin{equation}\label{SOCP-mu}\tag{\text{SOCP$_\mu$}}
	\begin{aligned}
	\min_{x,y}\ & F(x,y)+\mu\big(f(x,y)-V(x)\big) \quad
	\mathrm{s.t.\,}
	&\ (x,y)\in \Sigma_{\mathrm{SOC}},\ G(x, y) \leq 0.
	\end{aligned}	
	\end{equation}
When both the KKT condition and a certain second-order optimality condition hold for $y\in S(x)$, one has
\begin{equation} {\rm gph}\,S:=\big\{(x, y) \in \mathbb{R}^{n+m}: y\in S(x) \big\}\subseteq  \Sigma_{\mathrm{SOC}}\subseteq \Sigma_{\mathrm{KKT}}.\label{relations}
\end{equation}
In general, the inclusions above are strict. If the second inclusion is strict, i.e., the set $\Sigma_{\mathrm{KKT}}$ is strictly larger than the set $\Sigma_{\mathrm{SOC}}$, then obviously it is easier for a local optimal solution of
\eqref{BLPP} to be a solution to   \eqref{SOCP-mu} than to \eqref{KKTCP-mu}. This means that the partial calmness for the combined program with second-order optimality conditions is more likely to hold than the one for the combined program with first-order optimality conditions.

{For the  bilevel programming problem where the lower level is unconstrained, }
when we add the second-order optimality condition,
the partially penalized problem becomes a nonlinear semidefinite programming problem.
  For the general \eqref{BLPP} where the lower level problem is a constrained optimization problem,  there are several different second-order optimality conditions. We propose the corresponding combined program with each second-order optimality condition. Similar to the KKT approach where one minimizes over the original variables and the multipliers, we also propose some relaxed version of these second-order combined programs where multipliers are used as variables.

  Another difficulty of the value function or the combined approach is that the value function is usually nonsmooth and implicit. Since the set of second-order stationary points $\Sigma_{\mathrm{SOC}}$ is in general smaller than the set of first-order stationary points $\Sigma_{\mathrm{KKT}}$, it is more likely that the set of second-order stationary points coincides with $ {\rm gph}\,S$. In particular, if it happens that $\Sigma_{\mathrm{SOC}}= {\rm gph}\,S$, then the value function constraint $f(x,y)-V (x)\leq0$ can be removed from \eqref{SOCP} and so the partial calmness of the problem \eqref{SOCP} holds with penalty parameter $\mu=0$. Consequently, the resulting necessary optimality condition is much easier to obtain and does not involve the value function.
This is another advantage of using the combined program with second-order optimality conditions.

\textbf{Outline.}  The remaining part of the paper is organized as follows. In Section \ref{sec2}, we gather some preliminaries and preliminary results that will be used later.  An illustrative example will be given in Section \ref{example}. In Section \ref{sec4}, we introduce the combined problems with different kinds of second-order optimality conditions and the relaxed problems, discuss the partial calmness conditions and optimality conditions, and also give some examples. {Conclusions are given in Section \ref{sec5}.}

\textbf{Symbols and Notations.} Our notation is basically standard.
For a matrix $A$, we denote by $A^T$ its transpose. The inner product of two vectors $x,y$ is denoted by $x^T y$ or $\langle x,y\rangle$. We denote by $\mathbb{S}^{m}$ the set of symmetric $m \times m$ matrices equipped with the inner product $\langle A, B\rangle:=$ $\operatorname{tr}(A B)$, $A, B \in \mathbb{S}^{m},$ where $\operatorname{tr}(A)$ denotes the trace of the matrix $A$. The notation $A \succeq 0$ ($A \preceq 0$) means that $A$ is a symmetric positive (negative) semidefinite matrix. The set of symmetric positive semidefinite $m \times m$ matrices is denoted by $\mathbb{S}_{+}^{m}$.
For $z\in\mathbb{R}^d$ and  $\Omega\subseteq\mathbb{R}^d$, we denote by $\mathrm{dist}(z;\Omega)$ the distance from $z$ to $\Omega$.
For a { differentiable mapping}  $h:\mathbb{R}^d\rightarrow\mathbb{R}^{l}$, we denote its  Jacobian matrix   by $D h(z)\in \mathbb{R}^{l \times d}$ and by $\nabla h(z):=D h(z)^T $ the transpose of the Jacobian,  and in  the case where $l=1$,  its gradient vector by $\nabla h(z)$ and its Hessian   by $\nabla^2 h(z)$. For a nonsmooth function $g:\mathbb{R}^d\rightarrow\mathbb{R}$, we denote the {Clarke generalized gradient} of $g$ at $z$ by $\partial^c g(z)$. { For a set-valued mapping $\Gamma : K_1 \rightrightarrows K_2$ where $K_1$ and $K_2$ are Euclidean spaces, we denote the domain of $\Gamma$ with $\operatorname{dom} \Gamma :=\{x\in K_1 | \Gamma(x)\not= \emptyset\}$.}

\section{Preliminaries and preliminary results}\label{sec2}
In this section, we review and obtain some results that are needed in this paper.
\subsection{Second-order optimality conditions for the lower level program}

In this subsection, we review some results on second-order optimality conditions for the lower level program of \eqref{BLPP}.

{ We first recall some second-order optimality conditions for the following nonlinear optimization problem:
\begin{equation}\label{NP}\tag{NLP}
\begin{aligned}
\min_{t}&\  f(t)\ \mathrm{s.t.}&\ g(t) \leq 0, \ h(t)=0,
\end{aligned}
\end{equation}
where $f: \mathbb{R}^{m} \rightarrow \mathbb{R}$, $g: \mathbb{R}^{m} \rightarrow \mathbb{R}^{p_1}$, $h: \mathbb{R}^{m} \rightarrow \mathbb{R}^{p_2}$ are twice continuously differentiable.
Denote the Lagrangian function for \eqref{NP} by
\begin{equation*}
L(t,u,v):=f(t)+\sum_{i=1}^{p_1} u_{i} g_{i}(t)+\sum_{i=1}^{p_2} v_{i} h_{i}(t),\ \mathrm{for}\  (t,u,v)\in \mathbb{R}^{m} \times \mathbb{R}_{+}^{p_1}\times \mathbb{R}^{p_2},
\end{equation*}
and the generalized Lagrangian function for \eqref{NP} by
\begin{equation*}
\mathcal{L}_0(t, u_0, u,v):=u_{0} f(t)+\sum_{i=1}^{p_1} u_{i} g_{i}(t)+\sum_{i=1}^{p_2} v_{i} h_{i}(t),\ \mathrm{for}\  (t, u_{0}, u, v)\in \mathbb{R}^{m} \times \mathbb{R}_{+} \times \mathbb{R}_{+}^{p_1}\times \mathbb{R}^{p_2}.
\end{equation*}
Given a feasible point $t$ of problem \eqref{NP}, we denote the set of KKT multipliers at $t$ as follows:
\begin{equation*}
M^{1}(t):=\left\{(u,v)\in \mathbb{R}^{p_1}\times \mathbb{R}^{p_2}: \nabla_t L(t,u,v)=0, \ 	u \geq 0, \ \sum_{i=1}^{p_1} u_{i} g_{i}(t)=0 \right\}.
\end{equation*}
We define the critical cone at $t$ as follows:
\begin{equation*}
\mathcal{C}(t):=\Big\{d \in \mathbb{R}^{m}: D f(t) d \leq 0,\  D g_{j}(t) d \leq 0, \ \forall\,j \in J_{0}(t), D h_{i}(t) d = 0, i=1,...,p_2\Big\},
\end{equation*}
where $J_{0}(t):=\left\{j: g_{j}(t)=0\right\}$ denotes the set of indices of active inequalities at $t$. When $(u,v) \in M^{1}(t)$, by using the KKT condition, the critical cone can be written as
\begin{equation}
\label{KKTcritical}
\mathcal{C}(t)=\left\{d: \begin{array}{l}
	D g_{j}(t) d=0 \text { if } u_{j}>0,\  D g_{j}(t) d \leq 0 \text { if } u_{j}=0,\  \forall\,j \in J_{0}(t),\\
	 D h_{i}(t) d = 0, i=1,...,p_2 \end{array} \right\}.
\end{equation}
Another important set is the critical subspace given by
\begin{equation}\label{criticalsubspace}
\mathcal{S}(t):=\Big\{d \in \mathbb{R}^{m}: D g_{j}(t) d=0,\  \forall\,j \in J_{0}(t), D h_{i}(t) d = 0, i=1,...,p_2\Big\}.
\end{equation}
Note that when $M^{1}(t) \neq \emptyset$, the critical subspace $\mathcal{S}(t)$ is the lineality space of the critical cone $\mathcal{C}(t)$ (i.e., the largest linear space contained in $\mathcal{C}(t)$) and then $\mathcal{S}(t)=\mathcal{C}(t) \cap (-\mathcal{C}(t))$. If the strict complementarity holds, i.e., $u_j>0$, $\forall\,j\in J_0(t)$, we have $\mathcal{S}(t)=\mathcal{C}(t)$.

Now we review some classical second-order conditions.

\begin{definition}\label{DF-SOC}
	Let $t$ be a feasible point of problem \eqref{NP}. If $M^{1}(t) \neq \emptyset$, we say that
	\begin{enumerate}[(i)]
		\item the basic second-order optimality condition $(\mathrm{BSOC})$ holds at $t,$ if $\forall\,d \in \mathcal{C}(t)$, there exists $(u,v) \in M^{1}(t)$ such that $d^{T} \nabla_{tt}^{2} L(t,u,v) d \geq 0$;
		
		\item the weak second-order optimality condition $(\mathrm{WSOC})$ holds at $t$, if there exists $(u,v) \in M^{1}(t) \text { such that } d^{T} \nabla_{tt}^{2} L(t,u,v) d \geq 0, \forall\, d \in \mathcal{S}(t)$;
		
		\item the strong second-order optimality condition $(\mathrm{SSOC})$ holds at $t,$ if there exists $(u,v) \in M^{1}(t)$ such that $d^{T} \nabla_{tt}^{2} L(t,u,v) d \geq 0, \forall\, d \in \mathcal{C}(t)$.
\end{enumerate}
\end{definition}
}

  Note that when the linear independence constraint qualification (LICQ) holds at a feasible point $t$, there is a unique multiplier, i.e., the set $M^1(t)$ is a singleton. Hence, BSOC is equivalent to SSOC under LICQ.
All KKT type second-order optimality conditions such as BSOC, WSOC and SSOC hold at (local) minimizers only if certain constraint qualifications  are valid. BSOC requires a fairly weak constraint qualification. In classical results, MFCQ was required for BSOC to hold, c.f.,  \cite[Proposition 5.48]{BS00}. Recently
under a much weaker constraint qualification called the directional metrical subregularity condition \cite[Theorem 5.2]{G13}, it was shown that BSOC holds.
However,  WSOC and SSOC require much stronger constraint qualifications. 
It is known that SSOC (and hence WSOC) holds under the relaxed constant-rank constraint qualification (RCRCQ) \cite[Theorem 6]{MS11} and it is known that MFCQ was shown to be not enough for SSOC to hold {\cite[page 1350]{A91}}. Another condition called the critical regularity condition, which is not stronger than RCRCQ, is enough to give SSOC at local minimizers \cite[Theorem 2.1]{ML16}. Recently, it was shown that WSOC holds under MFCQ plus the weak constant rank property \cite[Theorem 3.1]{BHRV18}.

Even when no constraint qualification is assumed, a Fritz John second-order optimality condition (FJSOC) always holds at a local minimizer.

\begin{theorem}\cite[Proposition 5.48]{BS00}\label{thmFJ}
	Suppose $t$ is a local minimizer of \eqref{NP}. Then, for all $d \in \mathcal{C}(t)$, there is a Fritz John multiplier $(u_{0},u,v)$ such that
	\begin{equation*}
	d^{T} \nabla_{tt}^{2} \mathcal{L}_0(t, u_0, u,v) d \geq 0.
	\end{equation*}
\end{theorem}


{ Since we will use the above concepts for the lower level program of BLPPs frequently, we give the following notation. For fixed upper variable $x$ of \eqref{BLPP}, we denote the Lagrangian function and the generalized Lagrangian function for the lower level program by $L(y, u; x)$ and $\mathcal{L}_0(y, u_0, u; x)$, respectively. For any $y \in S(x)$, we denote the set of KKT multipliers for the lower level program \ref{Lx} at $y$ by $M^{1}(y;x)$. For any $u\in M^{1}(y;x)$, we call $(y,u)$ a KKT pair of program \ref{Lx}. We use $\mathcal{C}(y ; x)$ and $\mathcal{S}(y ; x)$ to denote the critical cone and the critical subspace at $y$ for fixed $x$, respectively.}

	{Since it is difficult to deal with the set of indices of active inequalities in the definition of the critical cone,} we introduce slack variables $z:=(z_1,\dots,z_p)^T\in\mathbb{R}^p$ for the lower level program, and obtain
	\begin{equation}\label{Lx2}\tag*{\text{$\widetilde{P}(x)$}}
	\min _{y,z} f(x, y) \quad {\rm s.t. } \quad g(x, y)+z^2 = 0.
	\end{equation}
	{Here, $z^2:=(z_1^2,\dots,z_p^2)^T$.} The above problem is equivalent to \ref{Lx} in the following sense. For fixed $x$, if $y^*$ is a global (local) optimal solution of \ref{Lx}, then there exists $z^*$ such that $(y^*,z^*)$ is a global (local) optimal solution of \ref{Lx2}. Conversely, if $(y^*,z^*)$ is a global (local) optimal solution of \ref{Lx2}, then $y^*$ is a global (local) optimal solution of \ref{Lx}.
	
	 Let $(y,z)$ be a feasible point of problem \ref{Lx2}. By  definition, we say that $u$ is a multiplier and $(y,z, u)$ is a KKT triple of problem \ref{Lx2} provided that
	\begin{equation*}
	\nabla_{(y,z)}L(y,z,u;x)=0,
	\end{equation*}
	where
	$
	L(y,z,u;x):=f(x,y)+\sum_{i=1}^p u_i\big[g_i(x,y)+z_i^2\big].
	$
	That is,
	\begin{align*}
	&\nabla_y f(x,y)+\sum_{i=1}^p u_i \nabla_y g_i(x,y)=0,\\
	&u_i z_i=0,\ g_i(x,y)+z_i^2=0,\ i=1,\dots,p.
	\end{align*}
	Note that, different from the KKT multipliers in $M^1(y;x)$, the multipliers $u_i$ above are not necessarily nonnegative.
	
	 {Since the problem \ref{Lx2} has only equality constraints, if the KKT condition holds, then the critical cone and the critical subspace  of problem  \ref{Lx2} are equal and  given by
	\begin{equation}
	\mathcal{C}(y, z ; x)	=\mathcal{S}(y, z ; x):=\Big\{(d, \nu) \in \mathbb{R}^{m}\times\mathbb{R}^{p}:    D_{y} g_{i}(x, y) d + 2z_i \nu_i=0,\,\forall\,i \Big\}. \label{criticalconeforsp}
	\end{equation}}
As an optimization problem with equality constraints, WSOC and SSOC for problem \ref{Lx2}  coincide and hence we call it SOC. Let $(y,z,u)$ be a KKT triple of problem \ref{Lx2}.  We say that SOC holds at $(y,z,u)$ if 
\begin{equation}\label{soc}
(d,\nu)^T \nabla_{(y,z)}^2 L(y,z,u;x)(d,\nu) \geq 0, \quad \forall\,(d,\nu)\in \mathcal{C}(y, z ; x).
\end{equation}
Note that
		\begin{equation*}
		\nabla^2_{(y,z)} L(y,z,u;x)=	
		\begin{pmatrix}
		\nabla_{yy}^2 L(y,u;x)& 0
		\vspace{2mm}\\
		0& 2 \text{diag}(u)
		\end{pmatrix},
		\end{equation*}
where $\text{diag}(u)$ denotes the $p\times p$ diagonal matrix with the elements of vector $u$ on the main diagonal. Thus
		\begin{equation}
		(d, \nu)^T\nabla^2_{(y,z)} L(y,z,u;x)(d, \nu)
		=d^{T} \nabla_{y y}^{2} L(y, u;x) d +2\sum_{i=1}^p u_i \nu_i^2.\label{secondorderc}
		\end{equation}

It is a simple matter to show that if  $(y^*,u)$ is a KKT pair of \ref{Lx}  then there exists $z^*$ such that  $(y^*,z^*,u)$ is a KKT triple  of \ref{Lx2}.
	Moreover suppose that $(y^*,u)$ satisfies WSOC for \ref{Lx}. Then
	$$d^{T} \nabla_{y y}^{2} L(y, u;x) d\geq 0 \quad \forall d \in \mathcal{S}(y ; x).$$
	By (\ref{criticalconeforsp}) and  (\ref{criticalsubspace}), we have
	\begin{align*}
	(d, \nu) \in \mathcal{S}(y, z ; x) \Longrightarrow d \in \mathcal{S}(y ; x).
	\end{align*}
	{Since $u\geq0$ for KKT pair $(y^*,u)$ of \ref{Lx}, by \eqref{secondorderc}, the following result is valid.}
	
	\begin{proposition}\label{equiv1}
		Let $(y^*,u)$  be a KKT pair of \ref{Lx}. Then there exists $z^*$ such that $(y^*,z^*,u)$ is a KKT triple of \ref{Lx2}. Furthermore, if $(y^*,u)$ satisfies $\mathrm{WSOC}$ for \ref{Lx}, then $(y^*,z^*,u)$ satisfies $\mathrm{SOC}$ \eqref{soc} for \ref{Lx2}.
	\end{proposition}	
%

	But the converse is not always true, that is, even if $(y^*,z^*,u)$ is a KKT triple of \ref{Lx2}, $(y^*,u)$ is not necessarily a KKT pair of \ref{Lx}. In fact, the condition $u\geq0$, concerning the sign of the multiplier, may not hold. For a counterexample, we refer the reader to \cite[Example 3.2]{FF17}. Under the second-order sufficient conditions and {some constraint qualification}, it has been proved that KKT points of the original \ref{Lx} and the reformulated \ref{Lx2} problems are essentially equivalent, cf. \cite[Proposition 3.6]{FF17}. {  Moreover in the final remarks of \cite{FF17}, the authors asked if there are other conditions which guarantee equivalence of the KKT points. In the next result,  we answer this question by showing  that the converse holds under the second-order necessary condition and hence improve the result of \cite[Proposition 3.6]{FF17}.}
	
	\begin{proposition}\label{equiv2}
		Let $(y^*,z^*,u^*)$ be a KKT triple of \ref{Lx2}. Assume that $(y^*,z^*,u^*)$ satisfies SOC \eqref{soc}.
		Then $u^*_i\geq0$ for all $i=1,\dots,p$. Hence $(y^*,u^*)$ is a KKT pair of \ref{Lx} satisfying $\mathrm{WSOC}$.
	\end{proposition}	
	\begin{proof}
		First, since $(y^*,z^*,u^*)$ is a KKT triple of \ref{Lx2}, we have $u^*_i z^*_i=0$ for all $i=1,\dots,p$. Thus
		$
		u^*_i g_i(x,y^*)=-u^*_i(z^*_i)^2=0,
		$
	which implies that $u^*_i=0$ if $z^*_i\neq0$ or equivalently $g_i(x,y^*)\neq0$.
		
		Now we consider the index $j$ such that $g_j(x,y^*)=0=z^*_j$. Let us prove that in this case $u^*_j\geq0$.
%
		Taking $d^*=0$, $\nu^*_i=0$ for $i\neq j$ and $\nu^*_j=1$, by the formula for $\mathcal{S}(y^*, z^* ; x)$ in (\ref{criticalconeforsp}),  we have $(d^*,\nu^*)\in \mathcal{S}(y^*, z^* ; x)$.
		By (\ref{secondorderc}), we have
		\begin{equation*}
		0\leq(d^*, \nu^*)^T\nabla^2_{(y,z)} L(y^*,z^*,u^*;x)(d^*, \nu^*)=2u^*_j,
		\end{equation*}
		which implies that $u^*_j\geq0$. Hence, we conclude that $(y^*,u^*)$ is a KKT pair of \ref{Lx}.
		
		Next we show that $(y^*,u^*)$ satisfies WSOC. For every $d\in\mathcal{S}(y^*;x)$, we have $D_y g_j(x,y^*) d=0$ for all $j\in J_0(y^*;x)$. For $i\notin J_0(y^*;x)$, i.e., $z^*_i\neq0$, we take $\nu_i=-D_y g_i(x,y^*) d/(2z^*_i)$. For all $j\in J_0(y^*;x)$, take $\nu_j=0$. Then it is obvious that $(d,\nu)\in\mathcal{S}(y^*, z^* ; x)$. Hence by  (\ref{secondorderc})
		\begin{align*}
		0\leq(d, \nu)^T\nabla^2_{(y,z)} L(y^*,z^*,u^*;x)(d, \nu)
		&=d^{T} \nabla_{y y}^{2} L(y^*, u^*;x) d +2\sum_{i=1}^pu^*_i \nu_i^2\\
		&=d^{T} \nabla_{y y}^{2} L(y^*, u^*;x) d
		\end{align*}
		since $u^*_i=0$ for all $i\notin J_0(y^*;x)$ and $\nu_j=0$ for all $j\in J_0(y^*;x)$. Therefore, $(y^*,u^*)$ satisfies WSOC.
	\end{proof}	
	


\subsection{Lipschitz continuity of the value function and the upper estimate of the Clarke subdifferential of the value function}

For convenience, we quote the original result obtained by Gauvin-Dubeau in \cite{Gauvin-Dubeau} below. For  results under weaker assumptions and sharper upper estimates, the reader is referred to
\cite[Corollary 4.8]{GLYZ14} and \cite[Proposition 2]{Ye20}. Note that under extra assumptions, the convex hull operation in the formula below can be removed {  and a tighter bound for the subdifferential can be obtained; see e.g. \cite[Proposition 1]{Ye20} for the case where the lower level program is linear, and \cite[Section 5]{MNP12} for the case where the solution map $S$ is $V$-inner semicontinuous at the point of interest.  Note that in the last case, the uniform boundedness of $Y$ assumption can be removed}.
 \begin{proposition}\cite[Theorem 5.3]{Gauvin-Dubeau}
 	Assume that the set-valued map $Y$ is uniformly bounded around $\bar x$, i.e., there exists a neighborhood $U(\bar x)$ of $\bar x$ such that  $\cup_{x\in U(\bar x)} Y(x)$ is bounded. Suppose that MFCQ holds at each $y\in S(\bar x)$. Then the value function $V$  is Lipschitz continuous near $\bar x$ and the Clarke subdifferential of $V$ at $\bar x$ has the following upper estimate:
$$ \partial^c V(\bar x) \subseteq co \{\nabla_x f(\bar x, y')+ \nabla_x g(\bar x,y')^Tu': y'\in S(\bar x), u'\in M^1(y'; \bar x)\},$$
where $co\,C$ denotes the convex hull of the set $C$.  \end{proposition}


 \subsection{Constraint qualifications and optimality conditions for the combined problem}\label{cq}
As discussed in the introduction,  there are various reformulations of  \eqref{BLPP}. To simplify the discussion, in this section we  consider the following general combined problem:
\begin{equation}\label{GCP}\tag{GCP}
\begin{aligned}
\min_{x,y, u, w}&\  F(x,y) \\
\mathrm{s.t.}&\  f(x, y)-V(x)\leq 0,\ g(x,y)\leq 0,\ u\geq0,\ u^Tg(x,y)=0, \\
&\  H_1(x, y, u):=\nabla_y f(x,y)+\nabla_y g(x,y)^Tu=0, \\
&  H_2(x, y, u, w)\in C,
\end{aligned}
\end{equation}
where $x \in \mathbb{R}^{n}$, $y \in \mathbb{R}^{m}$, $u \in \mathbb{R}^{p}$, $w \in \mathbb{R}^{l}$ and the mappings $F, f: \mathbb{R}^{n} \times \mathbb{R}^{m} \rightarrow \mathbb{R},\ g: \mathbb{R}^{n} \times \mathbb{R}^{m} \rightarrow \mathbb{R}^p,$  $H_2: \mathbb{R}^{n} \times \mathbb{R}^{m} \times \mathbb{R}^{p} \times \mathbb{R}^{l} \rightarrow K$ are continuously differentiable, $K$ is a Euclidean space, and $C$ is a nonempty convex subset of $K$. Here the constraint $H_2(x, y, u, w)\in C$ represents a constraint which comes from a second-order optimality condition. Note that for simplicity we have omitted the upper level constraint in this section.

We define the partial calmness for \eqref{GCP} as follows.
\begin{definition}[Partial calmness for \eqref{GCP}]\label{pcdef}
	Let $(\bar{x}, \bar{y}, \bar{u}, \bar{w})$ be a local solution of \eqref{GCP}. We say that \eqref{GCP} is partially calm at $(\bar{x}, \bar{y}, \bar{u}, \bar{w})$ if there exists $\mu\geq0$ such that $(\bar{x}, \bar{y}, \bar{u}, \bar{w})$ is a local solution of the following partially penalized problem:
\begin{equation}\label{uniformp}\tag{\text{GCP$_\mu$}}
\begin{aligned}
\min_{x,y, u, w}&\  F(x,y)+\mu(f(x, y)-V(x)) \\
\mathrm{s.t.}&\ g(x,y)\leq 0,\ u\geq0,\ u^Tg(x,y)=0, \ { H_1(x, y, u)=0, H_2(x,y,u,w) \in C}.
\end{aligned}
\end{equation}
\end{definition}

Since the second-order condition $H_2(x, y, u, w)\in C$ is redundant, the feasible region of \eqref{uniformp} must include in the feasible region of the corresponding partially penalized problem (CP$_\mu$), and thus the partial calmness condition for problem \eqref{GCP} is easier to hold than that for problem \eqref{CP}.

\begin{proposition}\label{pcc}
		Let $(\bar{x},\bar{y},\bar u)$ be a local optimal solution to problem \eqref{CP}. Suppose that \eqref{CP} is partially calm at $(\bar{x},\bar{y},\bar u)$ and there is $\bar w$ such that $(\bar{x},\bar{y},\bar u, \bar w)$  is a local optimal solution for \eqref{GCP}, then \eqref{GCP} is also partially calm at $(\bar{x},\bar{y},\bar u, \bar w)$.\end{proposition}

We now study the constraint qualification and optimality condition for problem (\ref{GCP}). If the value function is Lipschitz continuous { and $K=\mathbb{R}^s$}, then problem (\ref{GCP}) is an MPEC with Lipschitz continuous problem data.
Due to the value function constraint, the nonsmooth MFCQ fails to hold at any feasible solution of the above problem \cite[Proposition 3.2]{YZ95}.

Recall that in MPEC literature,  one usually defines a Mordukhovich (M-) or a Strong { (S-) stationarity condition} based on whether  the multipliers are taken from the limiting normal cone or the regular normal cone of the complementarity set respectively. Similar to \cite[Definition 4.2]{Ye11}, we define M-/S- { stationarity condition} based on the value function for \eqref{GCP}. Given a feasible vector $(\bar{x}, \bar{y}, \bar{u}, \bar{w})$ of problem \eqref{GCP}, we define the following index sets:
\begin{eqnarray*}
	\begin{aligned}
		&\ I_{g}=I_{g}(\bar{x}, \bar{y}, \bar{u}, \bar{w}) :=\left\{j: g_{j}(\bar{x}, \bar{y})=0, \bar{u}_{j}>0\right\}, \\
		&\ I_{u}=I_{u}(\bar{x}, \bar{y}, \bar{u}, \bar{w}) :=\left\{j: g_{j}(\bar{x}, \bar{y})<0, \bar{u}_{j}=0\right\}, \\
		&\ I_{0}=I_{0}(\bar{x}, \bar{y}, \bar{u}, \bar{w}) :=\left\{j: g_{j}(\bar{x}, \bar{y})=0, \bar{u}_{j}=0\right\} .
	\end{aligned}
\end{eqnarray*}

\begin{definition}[Stationary conditions for \eqref{GCP} based on the value function]\label{defsm}
	Let $(\bar{x}, \bar{y}, \bar{u}, \bar{w})$ be a feasible solution to \eqref{GCP}.
	\begin{enumerate}[(i)]
		\item We say that $(\bar{x}, \bar{y}, \bar{u}, \bar{w})$ is an M-stationary point based on the value function if there exist $\mu \geq 0$, $\lambda^{g} \in \mathbb{R}^{p}$, $\lambda^{u} \in \mathbb{R}^{p}$ and $  { \lambda_1^H \in  \mathbb{R}^{m},\lambda_2^H} \in  \mathbb{R}^{s}$ such that
		\begin{align}
		&\ 0 \in  \Big[\nabla F(\bar{x}, \bar{y})+  \mu (\nabla f(\bar{x}, \bar{y})-\partial^c V(\bar x) \times\{0\} )  +\nabla g(\bar{x}, \bar{y})^T\lambda^{g} \Big] \times \big\{(0,0)\big\} \label{sc1}\\
                 &\ \quad\quad -  (0,0, \lambda^{u},0)+ \nabla H_1(\bar{x}, \bar{y}, \bar{u})^T\lambda_1^H+ \nabla H_2(\bar{x}, \bar{y}, \bar{u}, \bar{w})^T\lambda_2^H,\nonumber\\
		&\ \lambda_{j}^{g}=0, \  \forall\,j \in I_{u}, \quad  \lambda_{j}^{u}=0, \  \forall\,j \in I_{g},\quad  { \lambda_2^H\in  {N_{C}} (H_2(\bar{x}, \bar{y}, \bar{u}, \bar{w}))}, \label{sc2}\\
		&\ \text{and\ either}\  \lambda_{j}^{g}>0,\ \lambda_{j}^{u}>0, \ or\ \lambda_{j}^{g}\lambda_{j}^{u}=0,\  \forall\,j \in I_{0},\nonumber
		\end{align}
		where $N_{C}$ denotes the normal cone to the convex set $C$.
		
		\item We say that $(\bar{x}, \bar{y}, \bar{u}, \bar{w})$ is an S-stationary point based on the value function if there exist $\mu \geq 0, \lambda^{g} \in \mathbb{R}^{p}, \lambda^{u} \in \mathbb{R}^{p}$ and $  { \lambda_1^H \in  \mathbb{R}^{m},\lambda_2^H} \in  \mathbb{R}^{s}$ such that \eqref{sc1}$-$\eqref{sc2} and the following condition hold:
		\begin{eqnarray*}
			\begin{aligned}
               \lambda_{j}^{g}\geq0,\ \lambda_{j}^{u}\geq0,\  \forall\,j \in I_{0}.
			\end{aligned}
		\end{eqnarray*}
	\end{enumerate}
\end{definition}

{
By Definition \ref{defsm}, an M-/S- stationary point of problem \eqref{CP} must correspond to an M-/S- stationary point of problem \eqref{GCP}, but the converse is not true since the multiplier $\lambda_2^H$ can be nonzero.
When the multiplier $\lambda_2^H$ is nonzero, the combined approach of using only the first-order condition fails but the one using the second-order condition is useful.}

To obtain M-stationary conditions, we reformulate problem \eqref{GCP} equivalently as the following optimization problem:
\begin{equation}\label{gcpr}
\begin{aligned}
\min_{x,y, u, w}&\  F(x,y) \\
\mathrm{s.t.}&\  f(x, y)-V(x)\leq 0,\ (g(x,y),-u) \in \Omega_\mathrm{CS}^p,\\
&\ {  H_1(x, y, u)=0,\ H_2(x, y, u, w)\in C,}
\end{aligned}
\end{equation}
where
$\Omega_{\mathrm{CS}}^p:=\big\{(a,b)\in\mathbb{R}^{p}\times \mathbb{R}^{p}:a\leq0,b\leq0,\langle a,b\rangle=0\big\}$ is the negative complementarity set.

Denote the set of feasible solutions for problem \eqref{gcpr} by $\mathcal{F}$ and the perturbed feasible map by
\begin{equation}\label{fmap}
\mathcal{F}(r_1, r_2, r_3, { P_1, P_2}):=  \left\{(x,y, u, w):
\begin{aligned}
&\ f(x, y)-V(x)+r_1\leq 0, \\
&\ ( g(x,y)-r_2,-u+r_3 )\in \Omega_{\mathrm{CS}}^p, \\
&\ {  H_1(x, y, u)+P_1=0}, \\
&\ { H_2(x, y, u, w)+P_2 \in C}
\end{aligned}
\right\}.
\end{equation}
We now define the Clarke calmness for problem \eqref{GCP} as the one for its equivalent reformulation \eqref{gcpr} as follows.
\begin{definition}\label{ccdef}
	(Clarke calmness for problem \eqref{GCP}). Let $(\bar{x}, \bar{y}, \bar{u}, \bar{w})$ be a local optimal solution of \eqref{GCP}. We say that \eqref{GCP} is Clarke calm at $(\bar{x}, \bar{y}, \bar{u}, \bar{w})$ if there exist $\epsilon>0$ and $\mu\geq0$ such that, for all $(r_1, r_2, r_3,  { P_1, P_2})$ in $B(0, \epsilon)$, for all $(x, y, u, w)\in B((\bar{x}, \bar{y}, \bar{u}, \bar{w}),\epsilon)\cap\mathcal{F}(r_1, r_2, r_3, { P_1, P_2})$, one has
$$F(x, y)-F(\bar{x}, \bar{y})+\mu\|(r_1, r_2, r_3, { P_1, P_2})\|\geq0,$$
where $B(z,\epsilon)$ denotes the open ball centered at $z$ with radius $\epsilon$.
\end{definition}
Similar to  Burke  \cite[Theorem 1.1]{B91}, the Clarke  calmness defined in Definition \ref{ccdef} is equivalent to the exact penalization, i.e.,  \eqref{GCP} is Clarke calm at $(\bar{x}, \bar{y}, \bar{u}, \bar{w})$ if and only if there exists  $\mu\geq0$ such that it is a local solution of the penalized problem:
{\begin{equation*}\label{gcprp}
\begin{aligned}
\min_{x,y, u, w}  F(x,y)+\mu &\ (f(x,y)-V(x)+{  \sum_{i=1}^m |H_1^i(x, y, u)|} \\
&\ +\mathrm{dist}(( g(x,y),-u);\Omega_{\mathrm{CS}}^p)+\mathrm{dist}({ H_2}(x, y, u, w);C)).
\end{aligned}
\end{equation*}
}
It is well-known that the calmness of the perturbed feasible map (\ref{fmap})                                                                                                                                     or equivalently the existence of a local error bound for the feasible region $\mathcal{F}$  is a sufficient condition for Clarke calmness; see e.g. \cite[Proposition 2.2]{DSY14}. Moreover many classical constraint qualifications can be used to guarantee  the Clarke calmness at a local minimizer; see e.g. \cite[Proposition 2.3]{DSY14}.


Similar to the equivalence between the Clarke calmness and exact penalization, it was pointed out in  \cite[Proposition 3.3]{YZ95} that there is an equivalence between the partial calmness and partial exact penalization.
The Clarke calmness condition is in general stronger than the partial calmness. The partial calmness condition plus a usual constraint qualification for the partially penalized problem implies the Clarke calmness condition \cite[Theorem 3.1]{YZ95}. One may derive sufficient condition for the calmness for the general combined program using the results on the relaxed constant positive linear dependence constraint qualification {(RCPLD) \cite[Theorem 3.2]{MM21}, \cite[Theorem 3.2]{XY20}.}

We can now state the optimality conditions for the general combined program below. In fact, one can also apply the directional calmness and optimality conditions in  \cite[Theorem 3.1]{BY20}, which was developed using the directional approach to variational analysis  in \cite{BGO19},  to the general combined problem. To obtain S-stationary condition, we introduce the following constraint qualification.
\begin{definition}[MPEC LICQ]
	Let $(\bar{x}, \bar{y}, \bar{u}, \bar{w})$ be a feasible solution to problem \eqref{uniformp}. We say that MPEC LICQ holds at  $(\bar{x}, \bar{y}, \bar{u}, \bar{w})$ if
	the following non-degeneracy  condition holds:	
	\begin{eqnarray*}
&& \left \{ \begin{array}{l}
 \begin{aligned}0 =&\ \displaystyle \sum_{j\in J_0(\bar{x}, \bar{y})} \lambda_j^g \nabla g_j(\bar x,\bar y) \times \{(0,0)\}{-\left\{\left(0,0,\sum_{j\in I_u\cup I_0}\lambda_j^ue_j, 0\right) \right\}} \\
 &\ { +\nabla H_1(\bar{x}, \bar{y}, \bar{u})^T\lambda_1^H \times \{0\}+ \nabla H_2(\bar{x}, \bar{y}, \bar{u}, \bar{w})^T\lambda_2^H},\end{aligned}\\
	  { \lambda_2^H \in \operatorname{span} N_{C}( H_2(\bar{x}, \bar{y}, \bar{u}, \bar{w})),}
	  \end{array} \right .\\
	 && \Rightarrow (\lambda^g,\lambda^u, { \lambda_1^H, \lambda_2^H)=(0,0,0,0)},
	 \end{eqnarray*}
where $e_j\in\mathbb{R}^p$ denotes the vector whose j-th component is 1, and others are all zero, and  $\operatorname{span}(\Pi)$ denotes the affine hull of the set $\Pi$.
\end{definition}

\begin{theorem}	\label{stationary}
Let $(\bar{x}, \bar{y}, \bar{u}, \bar{w})$ be a local optimal solution to \eqref{GCP}.
Suppose that \eqref{GCP}  is Clarke calm at $(\bar{x}, \bar{y}, \bar{u}, \bar{w})$, then $(\bar{x}, \bar{y}, \bar{u}, \bar{w})$ is an M-stationary point based on the value function. If \eqref{GCP} is partially calm at $(\bar{x}, \bar{y}, \bar{u}, \bar{w})$, either  $\mu=0$ or the value function is smooth,  and
MPEC LICQ
holds,
 then $(\bar{x}, \bar{y}, \bar{u}, \bar{w})$ is an S-stationary point based on the value function.
\end{theorem}
{
\begin{proof}

Since \eqref{GCP} is equivalent to \eqref{gcpr}, by \cite[Theorem 2.1]{DSY14} and  the expression for the limiting normal cone of the complementarity set, we get the result for the M-stationary point.  Similarly,   by Corollary 6 in \cite{GYZ20}  and  the expression for the regular normal cone of the complementarity set,  we get the result for the S-stationary point. Alternatively, if the set $C$ is a polyhedral set, then  we can also use the  \cite[Theorem 3.8]{Me20} to derive the desired result.
\end{proof}
\begin{remark} By definition, it is easy to see that if \eqref{GCP} is partially calm and problem \eqref{uniformp} is Clarke calm at $(\bar{x}, \bar{y}, \bar{u}, \bar{w})$, then the Clarke calmness for \eqref{GCP} holds.
\end{remark}

}


\section{An illustrative example}\label{example}

To illustrate the difficulties of BLPPs and our approach, we consider the following example for which all known approaches fail.

\begin{example}\label{unc}
	\begin{equation} \label{eg1}
	\begin{aligned}
	\min_{x,y}&\  \left(x-\frac{1}{2}\right)^2+y^2\quad
	\mathrm{s.t.}&\ y\in S(x):=\argmin\limits_{y}\left\{\frac{1}{4} y^{4}-\frac{1}{2} x y^{2}:y\in\mathbb{R}\right\}.
	\end{aligned}
	\end{equation}
\end{example}

\begin{figure}[!htbp]
	\begin{minipage}[t]{0.45\textwidth}
		\centering
		\includegraphics[height=3.5cm,width=5.5cm]{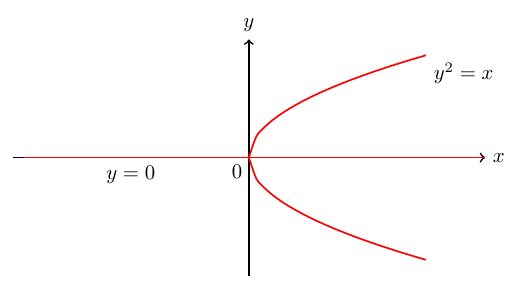}
		\scriptsize\caption{{Feasible set of problem \eqref{uncp}}}
		\label{fig1}
	\end{minipage}
	\hspace{.1in}
	\begin{minipage}[t]{0.45\textwidth}
		\centering
		\includegraphics[height=3.5cm,width=5.5cm]{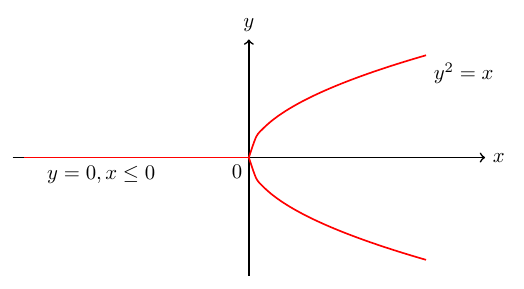}
		\scriptsize{\caption{Feasible set of problem \eqref{uncsp} (also the graph of $S(\cdot$))}}
		\label{figure2}
	\end{minipage}
\end{figure}

The first-order necessary condition for optimality of the lower level objective function with respect to $y$ is $y^3-xy=0$, which is equivalent to  saying that $y=0$ or $x=y^2$. Its graph is shown in Figure~\ref{fig1}.

Since the objective of the lower level program is not convex in { lower level variable} $y$, for each fixed $x$, not all corresponding $y$'s lying on the curve are global optimal solutions of the lower level program. The true global optimal solutions for the lower level problem are shown in Figure 2. It is easy to see that
\begin{eqnarray*}
	S(x)=\left\{\begin{array}{ll}
		\{\pm\sqrt{x}\} & {\rm if }\ x>0, \\
		\{0\} & {\rm if }\ x \leq 0,
	\end{array} \quad V(x)=\left\{\begin{array}{ll}
		-\frac{1}{4} x^{2} & {\rm if }\ x>0, \\
		0 & {\rm if }\  x \leq 0,
	\end{array}\right.\right.
\end{eqnarray*}
and $(\bar{x}, \bar{y})=(0,0)$ is the unique global optimal solution.

Now we claim that the partial calmness for \eqref{CP} does not hold at $(0,0)$. Indeed, the associated partially penalized problem  is given by
\begin{equation}\label{uncp}
\min _{x, y}\left\{F_{\mu}(x,y):=\left(x-\frac{1}{2}\right)^2+y^2+\mu\left(\frac{1}{4} y^{4}-\frac{1}{2} x y^{2}{-V(x)}\right): y^{3}-x y=0\right\}.
\end{equation}
Take any $\mu\geq 0$. For the objective values, we find
\begin{equation}
F_{\mu}\left(\frac{1}{k}, 0\right)=k^{-2}-k^{-1}+\frac{\mu}{4} k^{-2}+\frac{1}{4}, \quad F_{\mu}(0,0)=\frac{1}{4}.\label{pcfail}
\end{equation}
Thus, for $k>1+\mu / 4 ,\ F_{\mu}\left(\frac{1}{k}, 0\right)<F_{\mu}(0,0)$ holds, and this shows that  $(\bar{x}, \bar{y})=(0,0)$ is not a local minimizer of the associated partially penalized problem (\ref{uncp}). Hence the partial calmness for \eqref{CP} does not hold at $(0,0) .$ {  Moreover since \eqref{CP} is a standard nonlinear program, it is easy to check that the KKT condition does not hold at $(0,0)$.}


To explain our new approach, we now consider the following optimization problem in which we add the first and the second-order conditions to the value function reformulation of problem~\eqref{eg1}:
\begin{equation}\label{eg1v12}
\begin{aligned}
\min_{x, y}&\  \left(x-\frac{1}{2}\right)^2+y^2 \quad
\mathrm{s.t.}&\ f(x,y)-V(x)\leq0,\ y^{3}-x y=0,\ 3y^{2}-x \geq 0.
\end{aligned}
\end{equation}
Since both the first and the second-order conditions for the lower level program hold at $y\in S(x)$ without any further assumption, the constraints $y^{3}-x y=0$ and $3y^{2}-x \geq 0$ are redundant.
Hence $(\bar{x}, \bar{y})=(0,0)$ is still the optimal solution to the above problem.

From the graph in Figure 2, we can see that any point $(x,y)$ satisfying the first and the second-order conditions together lies in the graph of the solution mapping $S(\cdot)$. This means that the value function constraint can be removed and hence  $(0, 0)$ is a (local) minimizer of the following partially penalized problem with $\mu=0$:
\begin{equation}\label{uncsp}
\begin{aligned}
\min_{x, y}&\  \left(x-\frac{1}{2}\right)^2+y^2+\mu\big(f(x,y)-V(x)\big) \quad
\mathrm{s.t.}&\ y^{3}-x y=0,\ 3y^{2}-x \geq 0.
\end{aligned}
\end{equation}
Problem (\ref{uncsp}) is a one-level optimization problem. Furthermore, it is easy to check that its KKT condition holds at $(0,0)$.

Next we present a geometric explanation for Example \ref{unc}.


For Example~\ref{unc}, the partial calmness for \eqref{CP} at $(\bar{x}, \bar{y})=(0,0)$ means that for some $\mu\geq0$, $(\bar{x}, \bar{y})$ is still the optimal solution of the associated partially penalized problem \eqref{uncp}, whose feasible set is given by Figure~\ref{fig1}. But by (\ref{pcfail}), this is violated by taking points $\{(\frac{1}{k},0)\}_{k=1}^\infty$ on the line $\{(x,y):x>0, y=0\}$ in the feasible set.

To fix the above issue, we add the second-order necessary optimality condition of the lower level program in the combined problem \eqref{eg1v12}. The advantage of using the second-order necessary optimality condition is that the feasible set of the new associated partially penalized problem \eqref{uncsp} ruled out all of the points on the line $\{(x,y):x>0, y=0\}$ which are actually {local} maxima for the lower level objective function with $x>0$ (see Figure 2).

\section{Combined with second-order optimality conditions}\label{sec4}

A natural idea that comes from Example \ref{unc} is to add the second-order necessary optimality conditions of the lower level program in the combined problem. In this section, we consider combined problems with different kinds of second-order optimality conditions.

\subsection{Unconstrained case}

For the unconstrained bilevel programming problem
\begin{equation}\label{UBLPP}\tag{UBLPP}
\min _{x, y} F(x, y) \quad \text { s.t. } y \in \argmin _{y} f(x, y), \  G(x,y)\leq 0,
\end{equation}
we propose the following combined program using the second-order necessary optimality condition:
\begin{equation}\label{CPSOC}\tag{CPSOC}
\begin{aligned}
\min_{x,y}\ & F(x,y) \\
\mathrm{s.t.\ \,} &f(x,y)-V(x)\leq0, \ \nabla_y f(x,y) =0, \ \nabla_{y y}^{2} f(x, y) \in \mathbb{S}_{+}^{m},\
  G(x,y)\leq 0.
\end{aligned}
\end{equation}

We denote the corresponding partially penalized problem for \eqref{CPSOC} (as in Definition \ref{pcdef}) by (CPSOC$_\mu$). The problem (CPSOC$_\mu$) is a nonlinear semidefinite optimization problem. To derive an optimality condition for it, we may apply some constraint qualification, e.g., Robinson's constraint qualification (or a generalized MFCQ) of nonlinear semidefinite optimization problems.

\begin{theorem}\label{unconoc}
	Let $(\bar{x}, \bar{y})$ be a local optimal solution to \eqref{UBLPP}. Suppose that the partial calmness for \eqref{UBLPP} holds with either $\mu=0$ or with $\mu>0$ and the value function $V$ is Lipschitz continuous near $\bar{x}$. Then under some constraint qualification, there exist $\Omega \in \mathbb{S}_{+}^{m}$, $\mu \geq 0$, $\alpha \in \mathbb{R}^{m}$, and { $\beta \in \mathbb{R}^{q}_+$} such that
	\begin{eqnarray*}
		\begin{aligned}
			&\ 0 \in \nabla F(\bar{x}, \bar{y}) +\mu(\nabla f(\bar{x}, \bar{y})-\partial^c V(\bar{x}) \times\{0\})+\nabla (\nabla_{y} f)(\bar{x}, \bar{y})^T\alpha-{D \nabla_{y y}^{2} f(\bar{x}, \bar{y})^{*}} \Omega \\
			& \qquad \qquad + \nabla G(\bar x,\bar y)^T\beta,\\
			&\ \langle\nabla_{y y}^{2} f(\bar{x}, \bar{y}),\Omega\rangle=0,\ {  \beta^T G(\bar x,\bar y)=0},
		\end{aligned}
	\end{eqnarray*}
	where
	\begin{eqnarray*}
		D \nabla_{y y}^{2} f(\bar{x}, \bar{y})^{*} \Omega:=\left(\left\langle\frac{\partial}{\partial x_{1}} \nabla_{y y}^{2} f(\bar{x}, \bar{y}), \Omega\right\rangle, \ldots,\left\langle\frac{\partial}{\partial y_{m}} \nabla_{y y}^{2} f(\bar{x}, \bar{y}), \Omega\right\rangle\right)^{T}.
	\end{eqnarray*}
\end{theorem}

%
%

%

\subsection{Constrained case}

In the constrained case, as we reviewed in Section \ref{sec2}, there are four kinds of second-order optimality conditions:  FJSOC, BSOC, SSOC, and WSOC.


\subsubsection{Combined with the Fritz John second-order optimality condition}

We say that $y \in Y(x)$ is an FJSOC-point if for all $d\in\mathcal{C}(y; x)$, there exists $\left(u_{0}, u\right) \neq 0$ such that
\begin{eqnarray}\label{FJSOCpoint}
	\begin{aligned}
		&\ u_0 \nabla_{y} f(x, y)+\nabla_y g(x,y)^Tu=0, \\
		&\ g(x, y)\leq0,\  \left(u_{0}, u\right) \geq 0,\  \sum_{i=0}^pu_i=1,\ u^Tg(x, y)=0, \\
		&\ d^{T} \nabla_{y y}^{2} \mathcal{L}_0(y, u_0, u; x) d \geq 0.
	\end{aligned}
\end{eqnarray}
By Theorem~\ref{thmFJ}, if $y\in S(x)$ then $y$ is an FJSOC-point for \ref{Lx}.

Now we define
\begin{eqnarray*}
	\begin{aligned}
		&\  \Sigma_{\mathrm{FJSOC}}:=\Big\{(x, y) \in \mathbb{R}^{n+m}:{  \exists (u_0,u)\not =0 \mbox{ s.t. } (\ref{FJSOCpoint}) \mbox{ hold.}
		} \Big\},
	\end{aligned}
\end{eqnarray*}
and consider the following combined problem with FJSOC:
\begin{equation}\label{FJSOCP}\tag{FJSOCP}
\begin{aligned}
\min_{x,y}&\  F(x,y) \quad
\mathrm{s.t.}&\ f(x, y)-V(x) \leq 0, \ (x,y)\in \Sigma_{\mathrm{FJSOC}},\ G(x, y) \leq 0.
\end{aligned}
\end{equation}

Since it is not easy dealing with the set of indices of active inequalities in the critical cone, we propose to use the following set to relax the critical cone:
\begin{equation}\label{relaxedcone}
\{d\in \mathbb{R}^m: D_{y} f(x, y) d \leq 0,\  u_jD_{y} g_{j}(x, y) d \leq 0, \ \forall\,j=1,\dots,p\} \supseteq {\cal C}(y;x),
\end{equation}
{ where $(u_0,u)$ is an FJ-multiplier.} Under the strict complementarity, $``\supseteq"$ becomes $``="$ in the above relationship. Hence $y\in S(x)$ implies that there are $(u_0, u, d)$ such that the following relaxed FJ system holds:
\begin{eqnarray}\label{FJSOCpointnew}
	\begin{aligned}
		&\ u_0\nabla_{y} f(x, y)+\nabla_y g(x,y)^Tu=0, \\
		&\ g(x, y)\leq0,\  \left(u_{0}, u\right) \geq 0,\  \sum_{i=0}^pu_i=1,\ u^Tg(x, y)=0, \\
		&\ d^{T} \nabla_{y y}^{2} \mathcal{L}_0(y, u_0, u; x) d \geq 0,\\
		& \ D_{y} f(x, y) d \leq 0,\  u_jD_{y} g_{j}(x, y) d \leq 0, \ \forall\,j=1,\dots,p.
	\end{aligned}
\end{eqnarray}
Denote by
\begin{equation}\label{FJfunc}
K(x,y):=\{(u_0,u,d)\in\Xi(x,y)|u_0\nabla_{y} f(x, y)+\nabla_y g(x,y)^Tu=0\},
\end{equation}
where
\begin{equation}\label{FJfunc1}
\Xi(x,y):=  \left\{(u_0,u,d):
\begin{aligned}
&\ (u_0, u)\geq0,\ \sum_{i=0}^pu_i=1,\ u^Tg(x,y)=0, \\
&\ d^T\nabla_{y y}^{2}\big[u_{0} f(x,y)+\sum_{i=1}^{p} u_{i} g_{i}(x,y)\big] d\geq0, \\
&\ D_{y} f(x, y) d \leq 0,\  u_jD_{y} g_{j}(x, y) d \leq 0, \ \forall\,j=1,\dots,p
\end{aligned}
\right\}.
\end{equation}
Then problem \eqref{FJSOCP} can be reformulated as the following problem equivalently.
\begin{equation}\label{FJSOCP-2}\tag{FJSOCP-2}
\begin{aligned}
\min_{x,y}&\  F(x,y) \\
\mathrm{s.t.}&\ f(x, y)-V(x) \leq 0,\ g(x,y)\leq0,\ G(x, y) \leq 0, \\
&\ 0\in\bigcup_{(u_0,u,d)\in \Xi(x,y)} \left \{u_0\nabla_{y} f(x, y)+\nabla_y g(x,y)^Tu\right \}. \\
\end{aligned}
\end{equation}
Problem \eqref{FJSOCP-2} can be considered as an optimization problem with implicit variables as studied in recent paper \cite{BM21}.
Since
problem \eqref{FJSOCP-2} is still not practical to solve, we consider its explicit version---the relaxed combined problem with FJ second-order condition:
\begin{equation}\label{R-FJSOCP}\tag{R-FJSOCP}
\begin{aligned}
\min_{x,y,u_0,u,d}&\  F(x,y) \\
\mathrm{s.t.}&\ f(x, y)-V(x) \leq 0, \ u_0\nabla_{y} f(x, y)+\nabla_y g(x,y)^Tu=0, \ g(x,y)\leq0, \\
&\ (u_0, u)\geq0,\ \sum_{i=0}^pu_i=1,\ u^Tg(x,y)=0,\ d^T\nabla_{y y}^{2}\big[u_{0} f(x,y)+\sum_{i=1}^{p} u_{i} g_{i}(x,y)\big] d\geq0, \\
&\ D_{y} f(x, y) d \leq 0,\  u_jD_{y} g_{j}(x, y) d \leq 0, \ \forall\,j=1,\dots,p,\ G(x, y) \leq 0.
\end{aligned}
\end{equation}



Recall that a set-valued map $K$ is inner semicompact at $(\bar{x}, \bar{y})$ with respect to $\mathrm{dom}{K}$ if for each sequence $\left\{(x_{k},y_k)\right\}_{k \in \mathbb{N}} \subseteq \mathrm{dom}{K}$ such that $(x_{k},y_k) \rightarrow (\bar{x}, \bar{y})$, there is a convergent sequence $\left\{(u_0^\zeta,u^\zeta,d^\zeta)\right\}_{\zeta \in \mathbb{N}}$ and a subsequence $\left\{(x_{k}^{\zeta},y_{k}^\zeta)\right\}_{\zeta \in \mathbb{N}}$ such that $(u_0^\zeta,u^\zeta,d^\zeta) \in {K}\left(x_{k}^{\zeta},y_{k}^\zeta\right)$ holds for all $\zeta \in \mathbb{N}$; see e.g. \cite[page 7]{BM21}.   Note that the mapping $K$ defined in \eqref{FJfunc} is  automatically inner semicompact  since the set of FJ multipliers $(u_0,u)$ are uniformly bounded, and the sequence $d^\zeta$ in the definition can always be taken as zero. Hence by  \cite[Theorems 4.3 and 4.5]{BM21} we have the following equivalence between the problem with implicit variables and its explicit form.
\begin{proposition}\label{equi}
	Let $(\bar{x}, \bar{y})$ be a local (global) optimal solution to \eqref{BLPP}. Then for each $(u_0,u,d)\in K(\bar{x}, \bar{y}), (\bar{x}, \bar{y}, u_0, u, d)$ is a local (global) optimal solution of \eqref{R-FJSOCP}. Conversely, let $(\bar{x}, \bar{y}, u_0, u, d)$ be a global optimal solution to \eqref{R-FJSOCP}. Then $(\bar{x}, \bar{y})$ is a global solution of \eqref{BLPP}.  Moreover if for each $(u_0,u,d)\in K(\bar{x}, \bar{y})$, $(\bar{x}, \bar{y}, u_0, u, d)$ is a local optimal solution to \eqref{R-FJSOCP}, then $(\bar{x}, \bar{y})$ is a local solution of \eqref{BLPP}.
\end{proposition}


As in Definition \ref{pcdef}, we can define partial calmness for \eqref{FJSOCP} and partial calmness for \eqref{R-FJSOCP}, and denote the corresponding partially penalized problems by (FJSOCP$_\mu$) and (R-FJSOCP$_\mu$),
respectively.

Different from the relation between the partial calmness condition for \eqref{CP_{FJ}} and the partial calmness condition for \eqref{CPFJ}   in \cite[Theorem 4.4]{KYYZ}, the partial calmness condition for \eqref{FJSOCP} could not imply the partial calmness condition for \eqref{R-FJSOCP} directly because the critical cone has been relaxed in \eqref{R-FJSOCP}. But as we will show in Proposition \ref{propFJ}, the partial calmness condition for \eqref{CPFJ} implies the partial calmness condition for \eqref{R-FJSOCP}. {Note that this relation coincides with the result in Proposition \ref{pcc}.} On the other hand, since $\Sigma_{\mathrm{FJSOC}}\subseteq\Sigma_{\mathrm{FJ}}$, it is immediate that
\begin{eqnarray}\label{rlt}
\begin{aligned}
\text{partial calmness for \eqref{CP_{FJ}}}
\implies
\text{partial calmness for \eqref{FJSOCP}},
\end{aligned}
\end{eqnarray}
where $\Sigma_{\mathrm{FJ}}$ denotes the set of points which satisfy the Fritz John condition.

 In the following proposition, we show that the partial calmness for \eqref{R-FJSOCP} at $(\bar x, \bar y, \bar u_0, \bar u, \bar d)$ with $\bar d\neq0$ is { not stronger} than the one for \eqref{CPFJ} at $(\bar x, \bar y, \bar u_0, \bar u)$. Hence when the critical cone $\mathcal{C}(\bar y ; \bar x)\neq\{0\}$, one can always take a nonzero critical direction $\bar d$ to obtain a combined program with weaker partial calmness condition.
\begin{proposition}\label{propFJ}
Let $(\bar x,\bar y,\bar u_0,\bar u)$ be a local solution of \eqref{CPFJ}.
 Suppose that the partial calmness condition for \eqref{CPFJ} holds at $(\bar{x}, \bar{y}, \bar{u}_{0}, \bar{u})$. Then, for each $\bar{d}\in\mathcal{C}(\bar y ; \bar x)$, $(\bar{x}, \bar{y}, \bar{u}_{0}, \bar{u}, \bar{d})$ is a local optimal solution of problem  \eqref{R-FJSOCP} where the partial calmness condition holds. Conversely, suppose that problem \eqref{R-FJSOCP} is partially calm at a local solution $(\bar{x}, \bar{y}, \bar{u}_{0}, \bar{u}, 0)$ and $(\bar{x}, \bar{y}, \bar{u}_{0}, \bar{u})$ is a local solution of problem \eqref{CPFJ}. Then problem \eqref{CPFJ} is partially calm at $(\bar{x}, \bar{y}, \bar{u}_{0}, \bar{u})$.
	
\end{proposition}

\begin{proof}
{
Adding an additional variable $d$ and a superfluous constraint $d\in\mathbb{R}^m$ to problem \eqref{CPFJ}, the first assertion follows from Proposition \ref{pcc}.
}
%

Now suppose that problem \eqref{R-FJSOCP} is partially calm at $(\bar{x}, \bar{y}, \bar{u}_{0}, \bar{u}, 0)$. Then there exist $\mu \geq 0$ and a neighborhood $U(\bar x,\bar y, \bar u_0,\bar u, 0)$ of $(\bar{x}, \bar{y}, \bar{u}_{0}, \bar{u},0)$ such that
$$ F(\bar x,\bar y)\leq F(x,y)+\mu (f(x,y)-V(x)), \quad \forall (x,y,u_0,u,d)\in {\cal F}_R\cap U(\bar x,\bar y, \bar u_0,\bar u,0),$$
where ${\cal F}_R$ is the feasible region of problem (R-FJSOCP$_\mu$). Let $(x,y, u_{0}, {u})\in U(\bar x,\bar y, \bar u_0,\bar u)$ be a feasible solution of problem \eqref{CP2}. Then $(x,y,u_0,u,0)$ is feasible to problem (R-FJSOCP$_\mu$). Hence it follows that the problem \eqref{CPFJ} is partially calm at $(\bar{x}, \bar{y}, \bar{u}_{0}, \bar{u})$.
\end{proof}

{   Note that \eqref{R-FJSOCP} and its partially penalized problem (R-FJSOCP$_\mu$) are MPECs. Based on the discussions on the partial calmness condition,  Theorem \ref{stationary} can be used to derive S-/M- type necessary optimality conditions for problem \eqref{R-FJSOCP} under appropriate constraint qualifications.}

\subsubsection{Combined with the basic second-order optimality condition}
As reviewed in Section \ref{sec2}, under certain constraint qualifications, $M^1(y;x)\not =\emptyset$ for $y\in S(x)$ and one of the second-order optimality conditions BSOC, WSOC, and SSOC holds. In this subsection, we study the combined problem with BSOC. We say that $y$ is a BSOC, WSOC, or SSOC point of \ref{Lx} respectively if Definition \ref{DF-SOC}(i),  \ref{DF-SOC}(ii), or  \ref{DF-SOC}(iii) holds respectively.
%
%
Now we define
\begin{eqnarray*}
	\begin{aligned}
		&\Sigma_{\mathrm{BSOC}}:=\Big\{(x, y) \in \mathbb{R}^{n+m}:y\text{ is a BSOC-point for }P(x) \Big\},\\
		&\Sigma_{\mathrm{WSOC}}:=\Big\{(x, y) \in \mathbb{R}^{n+m}:y\text{ is a WSOC-point for }P(x) \Big\},\\
		&\Sigma_{\mathrm{SSOC}}:=\Big\{(x, y) \in \mathbb{R}^{n+m}:y\text{ is an SSOC-point for }P(x) \Big\}.
	\end{aligned}
\end{eqnarray*}
It is easily seen that
\begin{equation}
\Sigma_{\mathrm{SSOC}}\subseteq\Sigma_{\mathrm{BSOC}},\
\Sigma_{\mathrm{SSOC}}\subseteq\Sigma_{\mathrm{WSOC}},\
\mathrm{and}\ \Sigma_{\mathrm{SSOC}}\stackrel{\mathrm{LICQ}}{=}\Sigma_{\mathrm{BSOC}}.
\end{equation}

Similar to the combined problem \eqref{FJSOCP},  we consider the combined problem with basic (weak, strong) second-order optimality conditions \eqref{SOCP}
where $\Sigma_{\mathrm{SOC}}=\Sigma_{\mathrm{BSOC}}, \Sigma_{\mathrm{WSOC}}, \Sigma_{\mathrm{SSOC}}$, respectively. Different from FJSOC, none of BSOC, WSOC, and SSOC is necessary without extra constraint qualifications. Thus this reformulation requires that BSOC, WSOC, and SSOC hold at the optimal solution of the lower level program. At least it requires that the KKT conditions hold at the { optimal solutions} of the lower level program (i.e., $M^1(y;x)\neq\emptyset$).


Since it is difficult to express the set of indices of active inequalities directly in the combined problem \eqref{SOCP} with $\Sigma_{SOC}=\Sigma_{BSOC}$ such that it is still an optimization problem with equality and inequality constraints, we relax the critical cone \eqref{KKTcritical} as
\begin{eqnarray}
 \mathcal{C}(y ; x)&\subseteq &\Big\{d \in \mathbb{R}^{m}:  {D_y f(x,y)d\leq 0,}\ u_jD_{y} g_{j}(x, y) d \leq 0, \ \forall\,j=1,\dots,p \Big\}\nonumber \\
&=& \Big\{d \in \mathbb{R}^{m}:  u_jD_{y} g_{j}(x, y) d = 0, \ \forall\,j=1,\dots,p \Big\},\label{relaxedcone1} \label{relaxedcone2}
\end{eqnarray}
where { $u$ is a KKT multiplier} and (\ref{relaxedcone2}) follows from
 $$0\geq D_y f(x,y)d=-\sum_{j\in J_0(y;x)} u_j D_{y} g_{j}(x, y) d\geq 0.$$
Hence we propose to consider  the following relaxed problem for the combined problem \eqref{SOCP} with $\Sigma_{\mathrm{SOC}}=\Sigma_{\mathrm{BSOC}}$:
\begin{equation}\label{SOCCP}\tag{R-BSOCP}
\begin{aligned}
\min_{x,y,u,d}&\  F(x,y) \\
\mathrm{s.t.}&\ f(x, y)-V(x) \leq 0, \ \nabla_{y} f(x, y)+\nabla_y g(x,y)^Tu=0,\ g(x,y)\leq0, \\
&\ u\geq0,\ u^Tg(x,y)=0,\ d^T\nabla_{y y}^{2}\big[f(x,y)+\sum_{i=1}^{p} u_{i} g_{i}(x,y)\big] d\geq0, \\
&\ u_jD_{y} g_{j}(x, y) d = 0, \ \forall\,j=1,\dots,p,\ G(x, y) \leq 0.
\end{aligned}
\end{equation}


Similar to Proposition \ref{equi}, since there is the value function constraint, the combined problem \eqref{SOCP} and the relaxed combined problem \eqref{SOCCP} are both equivalent in global solutions to the original problem when the corresponding second-order optimality conditions hold \cite[Theorem 4.3]{BM21}. To state the relationship on local solutions, we define the following mapping:
\begin{equation}\label{Bfunc}
\widetilde{K}(x,y):=\{(u,d)\in\widetilde{\Xi}(x,y)|\nabla_{y} f(x, y)+\nabla_y g(x,y)^Tu=0\},
\end{equation}
where
\begin{equation}\label{Bfunc1}
\widetilde{\Xi}(x,y):=  \left\{(u,d):
\begin{aligned}
&\ u\geq0,\ u^Tg(x,y)=0, \\
&\ d^T\nabla_{y y}^{2}\big[f(x,y)+\sum_{i=1}^{p} u_{i} g_{i}(x,y)\big] d\geq0, \\
&\ u_jD_{y} g_{j}(x, y) d = 0, \ \forall\,j=1,\dots,p
\end{aligned}
\right\}.
\end{equation}

\begin{proposition}\label{equi2}
Let $(\bar{x}, \bar{y})$ be a local optimal solution to \eqref{BLPP}. Suppose that the basic second-order optimality condition holds for the lower level problem $P(\bar{x})$ at $\bar{y}$. Then for each $(u,d)\in \widetilde{K}(\bar{x}, \bar{y}), (\bar{x}, \bar{y}, u,d)$ is a local optimal solution of \eqref{SOCCP}. Conversely, let $(\bar{x}, \bar{y}, u,d)$ be a local optimal solution to \eqref{SOCCP} for each $(u,d)\in \widetilde{K}(\bar{x}, \bar{y})$, Furthermore, let $\widetilde{K}$ be inner semicompact at $(\bar{x}, \bar{y})$ with respect to $\mathrm{dom}\widetilde{K}$.
 Then $(\bar{x}, \bar{y})$ is a local solution of \eqref{BLPP}.
\end{proposition}

\begin{remark}
The function $\widetilde{K}$ defined in \eqref{Bfunc} is inner-semicompact when MFCQ for the lower level holds at $(\bar{x}, \bar{y})$. { This is true because at the optimal solutions, MFCQ is equivalent to the property of the set of KKT multipliers being nonempty and bounded \cite{Gau77}}.
\end{remark}

Next, we study the relation between the partial calmness for \eqref{SOCCP} and the partial calmness for \eqref{CP}. Similar to Proposition \ref{propFJ}, we can prove the following proposition.

\begin{proposition}\label{relationtoCP}
	Suppose that $(\bar x, \bar y, \bar u, 0)$ is a local solution of  \eqref{SOCCP}. Then the partial calmness for \eqref{SOCCP} holds at $(\bar x, \bar y, \bar u, 0)$ if and only if the partial calmness for \eqref{CP} holds at the local optimal solution $(\bar x, \bar y, \bar u)$.  Furthermore, if $\mathrm{BSOC}$ holds at $\bar y$ for the lower level problem $P(\bar x)$, then for all $\bar d\in\mathcal{C}(\bar y;\bar x)$, validity of the partial calmness condition for \eqref{CP} at the local optimal solution $(\bar x, \bar y, \bar u)$ implies validity of the partial calmness condition for \eqref{SOCCP} at the local optimal solution $(\bar x, \bar y, \bar u, \bar d)$.
\end{proposition}

{   Finally, \eqref{SOCCP} and its partially penalized problem (R-BSOCP$_\mu$) are MPECs.  Theorem \ref{stationary} can be used to derive S-/M- type necessary optimality conditions  under appropriate constraint qualifications.}

\subsubsection{Combined with the weak second-order optimality condition}\label{sec4.2.3}
{
If WSOC holds at the lower level, we can consider the following combined problem with WSOC:}
\begin{equation}\label{WSOCP}\tag{WSOCP}
\begin{aligned}
\min_{x,y,u}&\  F(x,y) \\
\mathrm{s.t.}&\ f(x, y)-V(x) \leq 0, \ \nabla_{y} f(x, y) +\nabla_y g(x,y)^Tu=0, \\
&\ g(x,y)\leq0,\ u\geq0,\ u^Tg(x,y)=0,\ G(x, y) \leq 0,\\
&\ 0\preceq\nabla_{y y}^{2}\big[f(x,y)+\sum_{i=1}^{p} u_{i} g_{i}(x,y)\big]\Big|_{\mathcal{S}(y ; x)}, \\
\end{aligned}
\end{equation}
and propose the corresponding partial calmness condition. Here
	$
	0\preceq\nabla_{y y}^{2}\big[f(x,y)+\sum_{i=1}^{p} u_{i} g_{i}(x,y)\big]\Big|_{\Gamma},
	$
with $\Gamma:=\mathcal{S}(y ; x)$ means that
$$d^T \nabla_{y y}^{2}\big[f(x,y)+\sum_{i=1}^{p} u_{i} g_{i}(x,y)\big] d\geq 0, \ \forall\,d\in \Gamma ,$$ i.e., the matrix $\nabla_{y y}^{2}\big[f(x,y)+\sum_{i=1}^{p} u_{i} g_{i}(x,y)\big]$  is a $\Gamma$-copositive matrix.

But the copositive matrix condition in \eqref{WSOCP} is not easy to tackle because the critical subspace ${\mathcal{S}(y ; x)}$ involves the set of indices of active inequalities of \ref{Lx}. To cope with this difficulty, the equivalence between the KKT points satisfying WSOC of the original problem \ref{Lx} and the reformulated problem \ref{Lx2} by introducing the squared slack variables is very useful. Indeed, by Propositions \ref{equiv1} and \ref{equiv2}, problem \eqref{WSOCP} is equivalent to the following reformulated problem by introducing the squared slack variables:
\begin{equation}\label{WSOCPZ}\tag{WSOCPZ}
\begin{aligned}
\min_{x,y,z,\lambda}&\  F(x,y) \\
\mathrm{s.t.}&\ f(x, y)-V(x) \leq 0, \ \nabla_{y} f(x, y)+\sum_{i=1}^p\lambda_i\nabla_{y} g_i(x, y)=0, \\
&\ g_i(x,y)+z_i^2=0,\ \lambda_i z_i=0,\ \forall\,i=1,\dots,p,\ G(x, y) \leq 0,\\
&\ 0\preceq\nabla^2_{(y,z)} L(y,z,\lambda;x)\Big|_{\mathcal{S}(y, z ; x)}. \\
\end{aligned}
\end{equation}
Now it is worth noting that the critical subspace
\begin{equation*}
\mathcal{S}(y, z ; x)=\Big\{(d, \nu) \in \mathbb{R}^{m}\times\mathbb{R}^{p}: \  D_{y} g_{i}(x, y) d + 2z_i \nu_i=0,\,\forall\,i \Big\}
\end{equation*}
does not involve the set of indices of active inequalities of \ref{Lx}.

Note that in contrast to all the other reformulations derived before, \eqref{WSOCPZ} is not a complementarity - but a generalized copositive programming problem with switching constraints \cite{LY21, Mehlitz20} which has slightly different properties than an MPEC.

\subsubsection{Combined with the strong second-order optimality condition}

 If SSOC holds at the lower level for each $y\in S(x)$, we can
consider the following combined problem:
	\begin{equation}\label{SSOCP}\tag{SSOCP}
	\begin{aligned}
	\min_{x,y,u}&\  F(x,y) \\
	\mathrm{s.t.}&\ f(x, y)-V(x) \leq 0,\ \nabla_{y} f(x, y)+\nabla_y g(x,y)^Tu=0, \\
	&\ g(x,y)\leq0,\ u\geq0,\ u^Tg(x,y)=0,\ G(x, y) \leq 0,\\
	&\ 0\preceq\nabla_{y y}^{2}\big[f(x,y)+\sum_{i=1}^{p} u_{i} g_{i}(x,y)\big]\Big|_{\mathcal{C}(y ; x)}. \\
	\end{aligned}
	\end{equation}
Recall that for a closed convex cone $\Gamma$, the class of all $\Gamma$-copositive matrices is the dual cone of the convex hull of $\big\{d d^T\in \mathbb{S}_{+}^{m}: d\in \Gamma\subseteq\mathbb{R}^m\big\}$ \cite[Lemma 2.28]{EJ08}. This provides a natural generalization of the constraint $\nabla_{y y}^{2} f(x, y) \in \mathbb{S}_{+}^{m}$ in the unconstrained case. {The problem \eqref{SSOCP} can be viewed as generalized semi-infinite programming problem \cite{S13,YW08} or generalized copositive programming problem (set-semidefinite optimization) \cite{BD12, EJ08}. { Similarly as in section \ref{sec4.2.3},  to deal with the difficulty of the active index set, one can also use the squared-slack-variable trick here. }}



\subsection{Examples and Summary}

In this section, we have discussed different types of combined problems with second-order optimality conditions, called \eqref{FJSOCP}, \eqref{SOCP}, \eqref{SSOCP} and \eqref{WSOCP}. To address the issue caused by the set of indices of active inequalities, we come up with the related relaxed problems, called \eqref{R-FJSOCP} and \eqref{SOCCP}, and also the problem with squared slack variables \eqref{WSOCPZ}. All of the combined and relaxed problems are equivalent to the original \eqref{BLPP} under some mild and necessary assumptions.

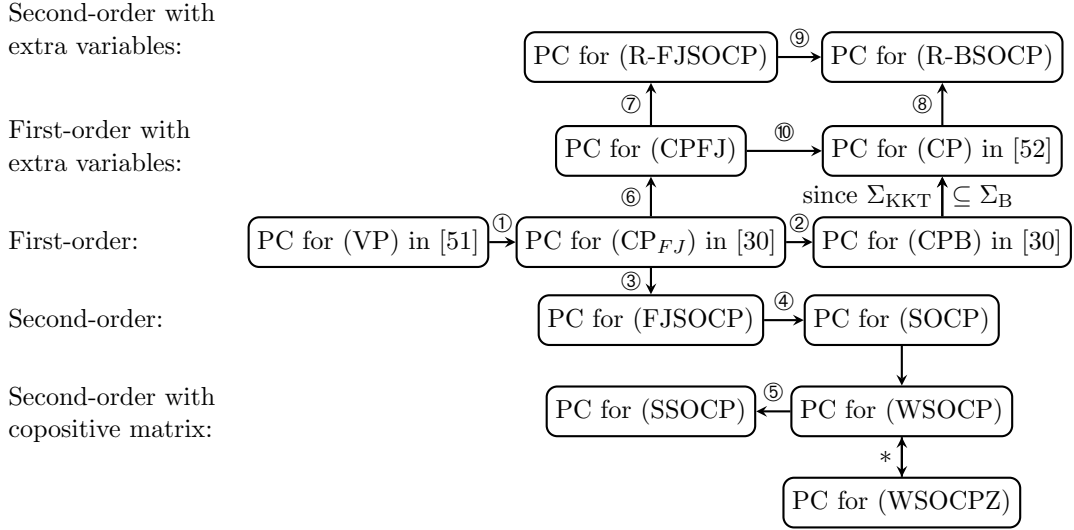
\begin{figure}[h!]
	\centering
\begin{tikzpicture}[thick,scale=0.6, every node/.style={scale=0.9}]
 \node[text width=3cm]                        (0)   {\makecell[l]{First-order:}};

 \node[draw, rounded corners, right=10pt of 0]                        (a)
 {\makecell[l]{PC for \eqref{VP} in \cite{YZ95}}};

 \node[draw, rounded corners, right=10pt of a]       (b)  {\makecell[l]{PC for \eqref{CP_{FJ}} in \cite{KYYZ} }};

 \node[draw, rounded corners, right=10pt of b]       (c)  {\makecell[l]{PC for \eqref{CPB} in \cite{KYYZ} }};


 \node[text width=3cm, above=10pt of 0]       (20)  {\makecell[l]{First-order with \\extra variables:}};

 \node[draw, rounded corners, above=15pt of c]       (d)  {\makecell[l]{PC for \eqref{CP} in \cite{YZ10}}};

 \node[draw, rounded corners, above=15pt of b]        (2b)  {\makecell[l]{PC for \eqref{CPFJ}}};

 \node[text width=3cm, below=10pt of 0]       (1)  {\makecell[l]{Second-order:}};

 \node[draw, rounded corners, below=10pt of b]        (3b)  {\makecell[l]{PC for \eqref{FJSOCP}}};

 \node[draw, rounded corners, right=15pt of 3b]        (3c)  {\makecell[l]{PC for \eqref{SOCP}}};

 \node[text width=3cm, above=10pt of 20]       (21)  {\makecell[l]{Second-order with\\ extra variables:}};

 \node[draw, rounded corners, above=16pt of 2b]        (4b)  {\makecell[l]{PC for \eqref{R-FJSOCP}}};

 \node[draw, rounded corners, above=16pt of d]        (4c)  {\makecell[l]{PC for \eqref{SOCCP}}};

  \node[text width=3cm, below=10pt of 1]       (22)  {\makecell[l]{Second-order with\\ copositive matrix:}};

  \node[draw, rounded corners, below=15pt of 3c]        (3d)  {\makecell[l]{PC for \eqref{WSOCP}}};

 \node[draw, rounded corners, below=15pt of 3b]        (3f)  {\makecell[l]{PC for \eqref{SSOCP}}};

  \node[draw, rounded corners, below=15pt of 3d]        (3e)  {\makecell[l]{PC for \eqref{WSOCPZ}}};

 \draw[arrow] (a)  -- node[anchor=south] {\ding{192}}  (b);
 \draw[arrow] (b)  -- node[anchor=south] {\ding{193}}  (c);
 \draw[arrow] (b)  -- node[anchor=east] {\ding{197}}  (2b);
 \draw[arrow] (2b)  -- node[anchor=south] {\ding{201}}  (d);
 \draw[arrow] (c)   -- node[left]  {since $\Sigma_{\mathrm{KKT}}$}  (d);
 \draw[arrow] (c)   -- node[right]  { $\subseteq\Sigma_{\mathrm{B}}$}  (d);
 \draw[arrow] (2b) -- node[anchor=east] {\ding{198}}  (4b);
 \draw[arrow] (d) -- node[anchor=east] {\ding{199}}  (4c);
 \draw[arrow] (4b) -- node[anchor=south] {\ding{200}}  (4c);
 \draw[arrow] (b) -- node[anchor=east] {\ding{194}}  (3b);
 \draw[arrow] (3b) -- node[anchor=south] {\ding{195}}  (3c);
 \draw[arrow] (3c) to  (3d);
  \draw[arrow] (3d) -- node[anchor=south] {\ding{196}}  (3f);
   \draw[arrow] (3d) to  (3e);
   \draw[arrow] (3e) -- node[anchor=east] {$*$}  (3d);


 \end{tikzpicture}
	\caption{Relationship between various partial calmness conditions. Here we denote ``partial calmness" briefly by PC. By Proposition \ref{pcc}, we have relations \ding{192}-\ding{196}. For relations \ding{197}, \ding{198}, and \ding{199}, we refer the reader to \cite[Theorem 4.4]{KYYZ}, Proposition \ref{propFJ}, and Proposition \ref{relationtoCP}, respectively. One may prove other relations by a similar argument of the proof of  Proposition \ref{propFJ}. The equivalent relation $*$ follows from Propositions \ref{equiv1} and \ref{equiv2}. An arrow between two PCs means one implies the other under certain constraint qualifications. Specifically, both the relations \ding{200} and \ding{201} require the validity of the KKT condition of \ref{Lx} for $y\in S(x)$.
}
	\label{relationship1}
\end{figure}

Similarly to \cite{KYYZ,YZ95,YZ10}, we have proposed various partial calmness conditions based on the combined problems above. We summarize the relationships between various partial calmness conditions in Figure \ref{relationship1}.

{   Under the validity of the  KKT condition, we can even establish the relationship between partial calmness for the FJ and the KKT type combined programs when the FJ multiplier considered satisfies $\bar u_0=0$. For example, even if the partial calmness for \eqref{CPFJ} holds for $(\bar x,\bar y,\bar u_0,\bar u)$ with $\bar u_0=0$, if the set of multiplier $M^1(\bar{y};\bar{x})$  is not empty and  $\widetilde{u}\in M^1(\bar{y};\bar{x})$, it can be shown that the partial calmness for \eqref{CP} holds for $(\bar x,\bar y, \widetilde{u}+k\bar u)$ when $k>0$ is sufficiently large. This comes from the fact that  $1+\sum_{j=1}^p(\widetilde{u}_j+k\bar u_j)=k+1+\sum_{j=1}^p\widetilde{u}_j$,  $(1,\widetilde{u}+k\bar u)/(1+\sum_{j=1}^p(\widetilde{u}_j+k\bar u_j))\rightarrow(0,\bar u)$ as $k\rightarrow+\infty$.}

Next, we use some nonconvex BLPPs to illustrate the combined approach with second-order optimality conditions and the necessary optimality conditions.

We first give an example for which the combined approach in \cite{KYYZ,YZ10} fails, but the partial calmness  {  and the necessary optimality condition
} will hold if one adds the basic second-order optimality condition for the lower level program in the associated combined problem.

\begin{example}\label{cons1}
	
	\begin{equation}\label{pcons1}
	\begin{aligned}
	\min_{x,y\in\mathbb{R}}&\  y^{2}-x\\
	\mathrm{s.t.}&\ -1 \leq x\leq 1,\ y\in S(x):=\argmin\limits_{y}\left\{\frac{1}{4} y^{4}-\frac{1}{2} x y^{2}:0 \leq y \leq \sqrt{2}\right\}.
	\end{aligned}
	\end{equation}

\medskip
\textbf{Claim:} In this example, we will show that
\begin{itemize}
\item{}  the partial calmness for \eqref{CP} does not hold at $(\bar{x}, \bar{y}, \bar{u})=(0, 0, 0)$;
\item{} the partial calmness for \eqref{SOCP} with $\Sigma_{\mathrm{SOC}}:=\Sigma_{\mathrm{BSOC}}=\Sigma_{\mathrm{SSOC}}$ holds at $(\bar{x}, \bar{y})=(0, 0)$;
\item{}  the partial calmness for \eqref{SOCCP} holds at $(\bar{x}, \bar{y}, \bar{u}, \bar{d})$ for any $\bar d\not =0$;
{ \item{} the partial calmness for \eqref{WSOCP} does not hold at $(\bar{x}, \bar{y}, \bar{u})=(0, 0, 0)$;}
{ \item{} necessary optimality conditions fail to hold for \eqref{CP};}
{ \item{} the S-stationarity condition holds for \eqref{SOCCP}.}
\end{itemize}
\medskip

	It is easy to see that
	\begin{eqnarray}\label{solution}
		S(x)=\left\{\begin{array}{ll}
		   \{\sqrt{2}\} & {\rm if }\ x>2, \\
			\big\{\sqrt{x}\big\} & {\rm if }\ 0<x\leq 2, \\
			\{0\} & {\rm if }\ x\leq  0,
		\end{array} \right. \quad
V(x)=\left\{\begin{array}{ll}
		    1-x & {\rm if }\ x > 2,\\
			-\frac{1}{4} x^{2} & {\rm if }\ 0<x\leq 2, \\
			0 & {\rm if }\ x\leq 0,
		\end{array}\right.
	\end{eqnarray}
	and $(\bar{x}, \bar{y})=(0,0)$ is a global optimal solution. Moreover, $M^1(0;0)=\{0\}$.

	Now we show that the partial calmness for \eqref{CP} does not hold at $(0,0,0)$. Indeed, the associated partially penalized problem  is given by
	\begin{equation}\label{2.1}
	\begin{aligned}
	\min_{x, y, u}&\ F_{\mu}(x,y):= y^{2}-x+\mu\left(\frac{1}{4} y^{4}-\frac{1}{2} x y^{2}-V(x)\right) \\
	\mathrm{s.t.}&\ y^{3}-xy-u_1+u_2=0,\ u_1\geq0,\ -u_1y=0, \\
	&\ u_2 \geq 0, u_2(y-\sqrt{2})=0,\ 0 \leq y \leq \sqrt{2},\ -1 \leq x\leq 1.
	\end{aligned}
	\end{equation}

	Note that when $x=\frac{1}{k} (k>0)$, $V(x)=-\frac{1}{4} x^{2}$. For any fixed $\mu$, the objective function value   $F_{\mu}(\frac{1}{k},0)=-k^{-1}+(\mu/4) k^{-2}<0=F_{\mu}(0,0)$ when $k>\mu / 4 .$ Hence $(\bar{x}, \bar{y}, \bar{u})=(0,0,0)$ is not a local minimizer of the associated partially penalized problem (\ref{2.1}) and the partial calmness for \eqref{CP} does not hold at $(0, 0,0)$.
	
Let us consider adding the second-order optimality conditions.  The critical cone is given by  \begin{eqnarray}\label{criticalcone}
		\mathcal{C}(y ; x)=\left\{\begin{array}{ll}
			\mathbb{R}_+ & {\rm if }\ x\in\mathbb{R},y =0, \\
			\mathbb{R} & {\rm if }\ x\in(0,2), y=\sqrt{x}, \\
			\{0\} & {\rm if }\ { x\geq2},y =\sqrt{2}.
		\end{array} \right.
	\end{eqnarray}
Since  LICQ holds, $\mathrm{BSOC}$ coincides with $\mathrm{SSOC}$ and hence $\Sigma_{\mathrm{BSOC}}=\Sigma_{\mathrm{SSOC}}$.  Problem \eqref{SOCP}	is given by
\begin{equation}\label{BSOCPex}
\begin{aligned}
\min_{x,y}&\  y^{2}-x \\
\mathrm{s.t.}&\ \frac{1}{4} y^{4}-\frac{1}{2} x y^{2}{-V(x)} \leq 0,\ (x,y)\in \Sigma_{\mathrm{SSOC}},\ -1 \leq x\leq 1.\\
\end{aligned}
\end{equation}
 Suppose $(x,y)\in \Sigma_{\mathrm{KKT}}$. Then it must satisfy the KKT condition
	\begin{eqnarray*}
	&& y^{3}-xy-u_1+u_2=0,\ u_1\geq0,\ -u_1y=0,\ u_2 \geq 0,\ u_2(y-\sqrt{2})=0,\ 0 \leq y \leq \sqrt{2},
	\end{eqnarray*}
	 with a unique multiplier $u$. It follows that
	\begin{equation}\label{KKT}
\Sigma_{\mathrm{KKT}}=\Big\{(x,0): x \in \mathbb{R}\Big\}\cup{ \Big\{(x,\sqrt{x}): x\in(0,2)\Big\}\cup\Big\{(x,\sqrt{2}): x \geq 2\Big\}}.
\end{equation}
	But $\mathrm{SSOC}$ states that
	$d^2(3y^2-x) \geq 0, \forall\,d \in \mathcal{C}(y ; x),$
which is equivalent to saying that
	$3y^2-x  \geq 0.$
	This means that the point $(x,y)$ with $x>0$ and $y=0$ does not satisfy $\mathrm{SSOC}$ and hence is not included in the set $\Sigma_{\mathrm{SSOC}}$. By the expression for the solution set \eqref{solution}, we have $\Sigma_{\mathrm{SSOC}}=\{(x,y): y\in S(x)\}$. Hence the value function constraint in problem \eqref{BSOCPex} holds for all $(x,y)\in \Sigma_{\mathrm{SSOC}}$. We therefore can remove the value function constraint from problem \eqref{BSOCPex}. This means that the partial calmness for \eqref{SOCP} with $\Sigma_{\mathrm{SOC}}=\Sigma_{\mathrm{SSOC}}$ holds at $(\bar{x}, \bar{y})=(0, 0)$ with $\mu=0$.


 Now consider the \eqref{SOCCP}:
	\begin{equation}\label{ex1rl}
	\begin{aligned}
	\min_{x, y, u, d}&\  y^{2}-x \\
	\mathrm{s.t.} & \frac{1}{4} y^{4}-\frac{1}{2} x y^{2}{-V(x)} \leq 0, y^{3}-xy-u_1+u_2=0,0 \leq y \leq \sqrt{2},\ -1 \leq x\leq 1,\\
    &\ u_1\geq0,\ -u_1y=0,\ -u_1d=0, u_2\geq0,\ u_2(y-\sqrt{2})=0,\ u_2d=0,\\
    &\ (3y^2-x) d^2 \geq 0.
	\end{aligned}
	\end{equation}
	Let $\bar d\not =0$. Then for any $d$ sufficiently close to $\bar d$, condition $(3y^2-x) d^2 \geq 0$ is equivalent to $3y^2-x \geq 0$. So similar to the analysis for the partial calmness for \eqref{SOCP}  with $\Sigma_{\mathrm{SOC}}=\Sigma_{\mathrm{SSOC}}$, the value function constraint can be removed.
	Then the partial calmness for  problem \eqref{SOCCP} holds at $(\bar{x}, \bar{y}, \bar{u}, \bar{d})=(0, 0, 0, \bar d )$ with $\mu=0$.


Recall that $\Sigma_{\mathrm{SSOC}}\subseteq \Sigma_{\mathrm{WSOC}}$ and the partial calmness with the larger set $\Sigma_{\mathrm{WSOC}}$ would be harder to hold. By the expression for the critical cone in (\ref{criticalcone}), we can obtain the expression for the critical subspace of the problem { \eqref{pcons1}}
	\begin{eqnarray*}
	\mathcal{S}(y ; x)=\left\{\begin{array}{ll}
		\{0\} & {\rm if }\ x\in\mathbb{R},y =0\ {\rm or}\ \sqrt{2}, \\
		\mathbb{R} & {\rm if }\ x\in\mathbb{R},0<y<\sqrt{2}.
	\end{array} \right.
	\end{eqnarray*}
	WSOC states that
	$$d^2(3y^2-x) \geq 0,\quad  \forall\,d \in \mathcal{S}(y ; x).$$
	Since when $x>0, y=0$, $d\in \mathcal{S}(y ; x)$ is taken as zero, these points are still in the set
	$\Sigma_{\mathrm{WSOC}}$ and hence $\Sigma_{\mathrm{WSOC}}=\Sigma_{\mathrm{KKT}}$. Since $\Sigma_{\mathrm{WSOC}}=\Sigma_{\mathrm{KKT}}$,
	for this example, the partial calmness for \eqref{SOCP} with $\Sigma_{\mathrm{SOC}}=\Sigma_{\mathrm{WSOC}}$ does not hold at $(\bar x,\bar y)$ and the partial calmness  for \eqref{WSOCP} does not hold at $(\bar x,\bar y,\bar u)$.	
	
	Point $(\bar{x}, \bar{y}, \bar{u})=(0, 0, 0)$ does not satisfy the stationary conditions for \eqref{CP} based on the value function as in Theorem \ref{stationary}.  Indeed, there do not exist $\mu\geq0$, $\beta$, $\eta^g$ and $\eta^G$ such that
	\begin{eqnarray*}
		\begin{aligned}
			0 \in&\  \nabla F(0,0)+\mu \left (\nabla f(0,0)-\partial^c V(0) \times\{0\} \right )+{ \nabla_{(x,y)} \left(\nabla _{y} L \right )(0,0;0)^T\beta} \\
			&\  + \nabla g(0,0)^T \eta^g+ \nabla G(0,0)^T \eta^G\\
		\end{aligned}
	\end{eqnarray*}
	since $\nabla F(0,0)=(-1,0)^T, \nabla  g_1(0,0)=(0,-1)^T$ and other terms are all zero.
	
	Problem \eqref{ex1rl} is an MPEC. The S-stationary condition based on the value function (Definition \ref{defsm}) holds at $(\bar{x}, \bar{y}, \bar{u}, \bar{d})=(0,0,0,1)$. Indeed, since $\nabla_{(x, y)}(\nabla_{y y}^{2} L)(0,0;0)\\=(-1,0)^T,$ there exists $\gamma=1$ (let other multipliers be all zero) such that
	\begin{eqnarray*}
		\begin{aligned}
			0 =&\  \nabla F(0,0)-{ \nabla_{(x, y)} \left(\nabla_{y y}^{2} L\right)(0,0;0)}^T \gamma.\\
		\end{aligned}
	\end{eqnarray*}

\end{example}





%
%

 In the following example, we show that the partial calmness for \eqref{WSOCP} may hold.
\begin{example}\label{cons2}
	\begin{equation}\label{pcons2}
	\begin{aligned}
	\min_{x,y}&\  y^{2}-x\\
	\mathrm{s.t.}&\ -1 \leq x\leq 1,\quad y\in S(x):=\argmin\limits_{y}\left\{\frac{1}{4} y^{4}-\frac{1}{2}xy^{2}:-1 \leq y \leq 1\right\}.
	\end{aligned}
	\end{equation}
{ 		It is easy to see that
	\begin{eqnarray}\label{}
		S(x)=\left\{\begin{array}{ll}
            \{\pm 1\} & {\rm if }\ x > 1,  \\
			\big\{\pm \sqrt {x}\big\} & {\rm if }\ 0<x\leq 1, \\
			\{0\} & {\rm if }\ x\leq 0,
		\end{array} \right.
	\end{eqnarray}
	\begin{eqnarray*}
\Sigma_{\mathrm{KKT}}=\Big\{(x,0): x \in \mathbb{R}\Big\}\cup{ \Big\{(x,\pm\sqrt{x}): x\in(0,1)\Big\}\cup\Big\{(x,\pm 1): x \geq 1\Big\}},
	\end{eqnarray*}
	and $(\bar{x}, \bar{y})=(0,0)$ is a global optimal solution. Moreover, $M^1(0;0)=\{0\}$.} Similarly to Examples \ref{unc} and \ref{cons1}, we can show that
\begin{itemize}
\item{} the partial calmness for \eqref{CP} {and necessary optimality conditions do} not hold at $(\bar{x}, \bar{y}, \bar{u})=(0, 0, 0)$;
\item{} the partial calmness for \eqref{SOCP} with $\Sigma_{\mathrm{SOC}}:=\Sigma_{\mathrm{BSOC}}=\Sigma_{\mathrm{SSOC}}$ holds at $(\bar{x}, \bar{y})=(0, 0)$;
\item{}  the partial calmness for \eqref{SOCCP} {and necessary optimality conditions hold} at $(\bar{x}, \bar{y}, \bar{u}, \bar{d})$ for any $\bar d\not =0$.
\end{itemize}
\medskip


%
%
	However, different from Example \ref{cons1}, we can show that the partial calmness for \eqref{WSOCP} also holds at $(\bar{x}, \bar{y}, \bar{u})=(0, 0, 0)$. In fact,  for problem  \eqref{pcons2},
{
	\begin{eqnarray*}
		\mathcal{S}(y ; x)=\left\{\begin{array}{ll}
			\{0\} & {\rm if }\ x\in\mathbb{R}, y =\pm1, \\
			\mathbb{R} & {\rm if }\ x\in\mathbb{R}, y\in (-1,1).
		\end{array} \right.
	\end{eqnarray*}}
But $\mathrm{WSOC}$ states that
	$d^2(3y^2-x) \geq 0,  \forall\,d \in \mathcal{S}(y ; x).$ Since
  points $(x,0)$ with $x>0$ do not satisfy the above $\mathrm{WSOC}$,
	we have $\Sigma_{\mathrm{WSOC}}={\rm gph}\,S$ (see  Figure 2).
	Hence the partial calmness for \eqref{WSOCP} holds at $(\bar{x}, \bar{y}, \bar{u})=(0, 0, 0)$ with $\mu=0$. 
%
	

%
\end{example}

%
%
%
%

We compare the results for the  two  examples in the following table.

\makeatletter\def\@captype{table}\makeatother
\begin{center}
	
	\scalebox{1}{
		\begin{tabular}{cccccc}
			\makecell[l]{\textbf{Examples}}
			& \makecell[l]{\textbf{\ref{CP}}}
			& \makecell[l]{\textbf{SOCP$_B$}}
            & \makecell[l]{\textbf{SOCP$_S$}}
            & \makecell[l]{\textbf{\ref{SOCCP}}}
			& \makecell[l]{\textbf{\ref{WSOCP}}}\\
			\toprule

			\makecell[l]{Example \ref{cons1}}
			&\makecell[l]{No}
			&\makecell[l]{Yes}
			& \makecell[l]{Yes}
            & \makecell[l]{Yes}
			& \makecell[l]{No}\\
			
			\midrule
			\makecell[l]{Example \ref{cons2}}
			&\makecell[l]{No}
			&\makecell[l]{Yes}
            & \makecell[l]{Yes}
			& \makecell[l]{Yes}
			& \makecell[l]{Yes}\\
			
			\bottomrule
	\end{tabular}}
	\caption{Comparison in the examples. Here we denote \eqref{SOCP} with with $\Sigma_{\mathrm{SOC}}=\Sigma_{\mathrm{BSOC}}$ or $\Sigma_{\mathrm{SSOC}}$  by \textbf{SOCP$_B$}, \textbf{SOCP$_S$}, respectively. ``Yes" or ``No" answers the question ``Does the partial calmness for the combined problem hold?"}
	\label{mainapps}
\end{center}

\section{Conclusions}\label{sec5}
In this paper,  we demonstrate that although the partial calmness condition is generic for a combined program with a first condition information,  there are still cases where the partial calmness condition {  and the corresponding necessary optimality conditions do} not hold.
To deal with these cases,  we propose to add both the first-order and  the second-order optimality conditions of the lower level problem as constraints.  There are several advantages in this approach. First,  by adding extra  constraints to the first-order combined problem, the new partial calmness condition and the resulting necessary optimality condition are easier to hold. Second,  by adding second-order optimality conditions,  it may be possible that the graph of the solution set to the lower level problem coincides with the set of second-order stationary points and hence the difficult value function constraint can be removed.  However there are also some drawbacks to our second-order approach.  First, by using a second-order optimality condition, we may need to introduce more extra variables other than the multipliers. Then similar to the difficulty in dealing with the reformulation involving with the KKT condition, the problem is no longer equivalent to the original bilevel program in the sense of local optimality.  Second,  since there are second-order optimality conditions in the feasible region of the partially penalized problem,  constraint qualifications may be harder to verify for the partially penalized problem.  Third, numerically the combined program with the second-order condition may sometimes be harder to solve than the combined program with the first-order condition. In our future work, we will try to utilize the advantages and  address the drawbacks.

\section*{Acknowledgments}{  We thank the guest editor Tim Hoheisel and the two anonymous referees for their suggestions, which have helped to improve the presentation of this paper.}

\end{document}